# Stability and Stabilization: Discontinuities and the Effect of Disturbances

Eduardo D. SONTAG*
*Department of Mathematics*
*Rutgers University*
*New Brunswick, NJ 08903, USA*
Email: `sontag@control.rutgers.edu`

## 1 Introduction

In this expository paper, we deal with several questions related to stability and stabilization of nonlinear finite-dimensional continuous-time systems. We review the basic problem of feedback stabilization, placing an emphasis upon relatively new areas of research which concern stability with respect to "noise" (such as errors introduced by actuators or sensors). The topics covered are as follows:

1. Review of Stability and Asymptotic Controllability.
2. The Problem of Stabilization.
3. Obstructions to Continuous Stabilization.
4. Control-Lyapunov Functions and Artstein's Theorem.
5. Discontinuous Feedback.
6. Nonsmooth CLF's.
7. Insensitivity to Small Measurement and Actuator Errors.
8. Effect of Large Disturbances: Input-to-State Stability.
9. Comments on Notions Related to ISS.

The discussion will be largely conceptual. Most technical details are left out, but links to precise statements in the literature are given. As a general reference, we direct the reader to the following Web site:

$$\boxed{\texttt{http:/www.math.rutgers.edu/~sontag}}$$

*The research was supported in part by Grant F49620-98-1-0242.

## 1.1   Classes of Systems Considered

We deal with rather general nonlinear finite-dimensional control systems of the following form:

$$\dot{x}(t) \;=\; f(x(t), u(t)) \tag{1}$$

where states $x(t) \in \mathbb{R}^n$ for all $t$ (the integer $n$ is the dimension of the system), and controls (or "inputs") are measurable locally essentially bounded maps

$$u(\cdot) \,:\, [0, \infty) \to \mathcal{U} = \mathbb{R}^m$$

into the control-value set, which we take here, for simplicity of presentation, to be simply an Euclidean space. (However, many of the results generalize easily to more general subsets $\mathcal{U} \subseteq \mathbb{R}^m$.) We assume that $f$ is locally Lipschitz on $(x, u)$, and $f(0,0) = 0$. In particular, we often consider systems with no inputs:

$$\dot{x}(t) \;=\; f(x(t))\,. \tag{2}$$

All definitions for such systems are implicitly applied as well to systems with inputs (1) by setting $u \equiv 0$; for instance, we define the global asymptotic stability (GAS) property for (2), but we say that (1) is GAS if $\dot{x} = f(x, 0)$ is.

The maximal solution $x(\cdot)$ of (1), corresponding to a given initial state $x(0) = x^0$, and to a given control $u$, is defined on some maximal interval $[0, t_{\max}(x^0, u))$ and is denoted by $x(t, x^0, u)$. For systems with no inputs (2) we write just $x(t, x^0)$. In general, we use the notation $|x|$ for Euclidean norms, and use $\|u\|$, or $\|u\|_\infty$ for emphasis, to indicate the essential supremum of a function $u(\cdot)$. For basic terminology and facts about control systems, we rely upon [36].

## 1.2   Stability and Asymptotic Controllability

The use of "comparison functions" has become widespread in stability analysis, as this formalism allows elegant formulations of most concepts. We recall the relevant definitions here. The class of $\mathcal{K}_\infty$ functions consists of all $\alpha : \mathbb{R}_{\geq 0} \to \mathbb{R}_{\geq 0}$ which are continuous, strictly increasing, unbounded, and satisfy $\alpha(0) = 0$, cf. Figure 1. The class of $\mathcal{KL}$ functions consists of those $\beta : \mathbb{R}_{\geq 0} \times \mathbb{R}_{\geq 0} \to \mathbb{R}_{\geq 0}$ with the properties that (1) $\beta(\cdot, t) \in \mathcal{K}_\infty$ for all $t$, and (2) $\beta(r, t)$ decreases to zero as $t \to \infty$. (It is worth remarking, cf. [38], that for each $\beta \in \mathcal{KL}$, there exist two functions $\alpha_1, \alpha_2 \in \mathcal{K}_\infty$ so that $\beta(r, t) \leq \alpha_2\left(\alpha_1(r) e^{-t}\right)$ for all $s, t$. Thus, as every function of the latter form is in $\mathcal{KL}$, and since $\mathcal{KL}$ functions are only used in order to express upper bounds, in a sense there is no need to introduce the class $\mathcal{KL}$.) We will also use $\mathcal{N}$ to denote the set of all nondecreasing functions $\sigma : \mathbb{R}_{\geq 0} \to \mathbb{R}_{\geq 0}$.

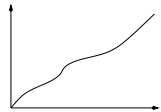

Figure 1: A function of class $\mathcal{K}_\infty$

Expressed in this language, the property of *global asymptotic stability (GAS)* for a system with no inputs (2) becomes:

$$(\exists\,\beta \in \mathcal{KL}) \quad |x(t,x^0)| \;\leq\; \beta\left(|x^0|,t\right) \quad \forall\, x^0,\; \forall\, t \geq 0\,.$$

It is an easy exercise to show that this definition is equivalent to the usual "$\varepsilon - \delta$" definition; for one implication, simply observe that

$$|x(t,x^0)| \;\leq\; \beta\left(|x^0|,0\right)$$

provides the stability (or "small overshoot") property, while

$$|x(t,x^0)| \;\leq\; \beta\left(|x^0|,t\right) \underset{t\to\infty}{\longrightarrow} 0$$

gives attractivity.

More generally, we define what it means for the system with inputs (1) to be (open loop, globally) *asymptotically controllable (AC)*. The definition amounts to requiring that for each initial state $x^0$ there exists some control $u = u_{x^0}(\cdot)$ defined on $[0,\infty)$, such that the corresponding solution $x(t,x^0,u)$ is defined for all $t \geq 0$, and converges to zero as $t \to \infty$, with "small" overshoot. Moreover, we wish to rule out the possibility that $u(t)$ becomes unbounded for $x$ near zero. The precise formulation is as follows.

$$(\exists\,\beta \in \mathcal{KL})\,(\exists\,\sigma \in \mathcal{N}) \quad \forall\, x^0 \in \mathbb{R}^n \;\; \exists\, u(\cdot),\; \|u\|_\infty \leq \sigma(|x^0|)\,,$$

$$|x(t,x^0,u)| \;\leq\; \beta\left(|x^0|,t\right) \quad \forall\, t \geq 0\,.$$

Finally, we say that $k: \mathbb{R}^n \to \mathcal{U}$ is a *feedback stabilizer* for the system with inputs (1) if $k$ is locally bounded (that is, $k$ is bounded on each bounded subset of $\mathbb{R}$), $k(0) = 0$, and the closed-loop system

$$\dot{x} \;=\; f(x, k(x)) \tag{3}$$

is GAS, i.e. there is some $\beta \in \mathcal{KL}$ so that $|x(t)| \leq \beta\left(|x(0)|,t\right)$ for all solutions and all $t \geq 0$. Obviously, if there exists a feedback stabilizer for (1), then (1) is also AC (just use $u(t) := k(x(t,x^0))$ as $u_{x^0}$). A most natural question is to ask if the converse holds as well, namely: *is every asymptotically controllable system also feedback stabilizable?* We will see that the answer is "yes", provided that we allow discontinuous feedbacks $k$, which in turn leads to the technical problem of defining precisely what one means by a "solution" of an initial value problem for (1) when $f(x, k(x))$ is not continuous (since, in that case, the standard theorems on existence do not apply). This is a question that is very much central to the rest of the paper.

## 1.3 Why Stabilization?

In control problems, there are often ultimate objectives other than stabilization. For example, one might merely be interested in making a *subset* of the state variables (as opposed to the full state) approach a desired set point, in which case (after translation of the origin of coordinates so as to make the desired values zero) one deals with partial or "output" stability: there is a given function $y = h(x)$ and one wishes to make $y(t) \to 0$ using suitable control actions. A variation of this idea is the formulation of tracking problems, where one wishes the output $y(t)$ to follow a predetermined path $y^*(t)$; this can be seen as an output stability problem for an "error" signal $e(t) = y(t) - y^*(t)$. Yet another variation consists in studying stability of invariant sets (as opposed to equilibria); for instance, one may be given a periodic orbit, and a feedback is sought so that for the closed-loop system it becomes a limit cycle. Finally, one may wish to allow for uncertainties in plant description (adaptive and/or robust control). Many such other objectives can be treated using the techniques that we shall describe; we focus on stabilization to an equilibrium because this is the simplest problem in which the main technical difficulties can be understood.

## 1.4 Regularity of Feedback

As already mentioned, one of the central issues with which we will be concerned is that of dealing with possibly discontinuous feedback laws $k$. Before addressing that subject, however, we will study what can be done with continuous feedback.

It turns out that requirements away from 0, say asking whether $k$ is continuous or smooth, are not very critical; it is often the case that one can "smooth out" a continuous feedback (or, even, make it real-analytic, via Grauert's Theorem) away from the origin. So, in order to avoid unnecessary complications in exposition due to nonuniqueness, let us call a feedback $k$ *regular* if it is locally Lipschitz on $\mathbb{R}^n \setminus \{0\}$. For such $k$, solutions of initial value problems $\dot{x} = f(x, k(x))$, $x(0) = x^0$, are well defined (at least for small time intervals $[0, \varepsilon)$) and, provided $k$ is a stabilizing feedback, are unique (cf. [36], Exercise 5.9.9).

On the other hand, behavior at the origin cannot be "smoothed out" and, at zero, the precise degree of smoothness plays a central role in the theory. For instance, consider the system

$$\dot{x} = x + u^3.$$

The continuous (and, in fact, smooth away from zero) feedback $u = k(x) := -\sqrt[3]{2x}$ globally stabilizes the system (the closed-loop system becomes $\dot{x} = -x$). However, there is no possible stabilizing feedback which is differentiable at the origin, since $u = k(x) = O(x)$ implies that

$$\dot{x} = x + O(x^3)$$

about $x = 0$, which means that the solution starting at any positive and small point moves to the right, instead of towards the origin. (A general result, assuming that $A$ has no purely imaginary eigenvalues, cf. [36], Section 5.8, is that if –and only if– $\dot{x} = Ax + Bu + o(x, u)$ can be locally asymptotically stabilized using a feedback which is differentiable at the origin, the linearization $\dot{x} = Ax + Bu$ must be itself stabilizable. In the example that we gave, this linearization is just $\dot{x} = x$, which is not stabilizable.)

## 2 Nonexistence of Regular Feedback

We now turn to the question of existence of regular feedback stabilizers. We first study a comparatively trivial case, namely systems with one state variable and one input. After that, we turn to multidimensional systems.

## 2.1 The Special Case $n = m = 1$

There are algebraic obstructions to the stabilization of $\dot{x} = f(x,u)$ if the input $u$ appears nonlinearly in $f$. Ignoring the requirement that there be a $\sigma \in \mathcal{N}$ so that controls can be picked with $\|u\| \leq \sigma(|x^0|)$, asymptotic controllability is, for $n = m = 1$, equivalent to:

$$(\forall x \neq 0)\,(\exists u)\ \ xf(x,u) < 0 \qquad (4)$$

(this is proved in [39]; it is fairly obvious, but some care must be taken to deal with the fact that one is allowing arbitrary measurable controls; the argument proceeds by first approximating such controls by piecewise constant ones). Let us introduce the following set:

$$\mathcal{O} := \{(x,u) \mid xf(x,u) < 0\},$$

and let $\pi : (x,u) \mapsto x$ be the projection into the first coordinate in $\mathbb{R}$. Then, (4) is equivalent to:

$$\pi\mathcal{O} = \mathbb{R} \setminus \{0\}.$$

(One can easily include the requirement "$\|u\| \leq \sigma(|x^0|)$" by asking that for each interval $[-K, K] \subset \mathbb{R}$ there must be some compact set $C_K \subset \mathbb{R}^2$ so that $[-K, K] \subseteq \pi(C_K)$. For simplicity, we ignore this technicality.) In these terms, a stabilizing feedback is nothing else than a locally bounded map $k : \mathbb{R} \to \mathbb{R}$ such that $k(0) = 0$ and so that $k$ is a section of $\pi$ on $\mathbb{R} \setminus \{0\}$:

$$(x, k(x)) \in \mathcal{O} \ \forall\, x \neq 0.$$

For a regular feedback, we ask that $k$ be locally Lipschitz on $\mathbb{R} \setminus \{0\}$.

Clearly, there is no reason for Lipschitz, or for that matter, just continuous, sections of $\pi$ to exist. As an illustration, take the system

$$\dot{x} = x\left[(u-1)^2 - (x-1)\right]\left[(u+1)^2 + (x-2)\right].$$

Let

$$\mathcal{O}_1 = \{(u+1)^2 < (2-x)\} \quad \text{and} \quad \mathcal{O}_2 = \{(u-1)^2 < (x-1)\}$$

(see Figure 2). Here, $\mathcal{O}$ has three connected components, namely $\mathcal{O}_2$ and $\mathcal{O}_1$ intersected

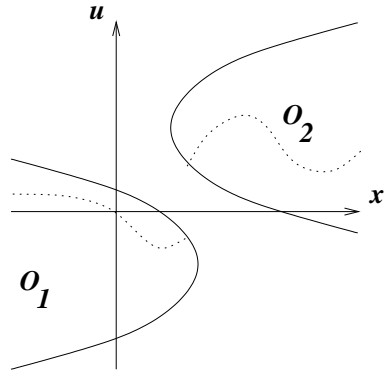

Figure 2: Two regions

with $x > 0$ and $x < 0$. It is clear that, even though $\pi\mathcal{O} = \mathbb{R}$, there is no continuous curve

(graph of $u = k(x)$) which is always in $\mathcal{O}$ and projects onto $\mathbb{R} \setminus \{0\}$. On the other hand, there exist many possible feedback stabilizers provided that we allow one discontinuity (e.g., dotted curve in figure). It is also possible to provide examples, even with $f(x, u)$ smooth, for which an infinite number of switches are needed in any possible stabilizing feedback law. The set $\mathcal{O}$ could be like the (non-shaded) set in Figure 3 (with an accumulation of smaller and smaller holes as we approach the origin). Finally, it may even be possible to stabilize semiglobally with a regular feedback, meaning that for each compact subset $K$ of the state-space there is a regular, even smooth, feedback law $u = k_K(x)$ such that all states in $K$ get driven asymptotically to the origin, but yet it may be impossible to find a single $u = k(x)$ which works globally. This illustrated by an example where the set $\mathcal{O}$ could be like the (non-shaded) set in Figure 4 (one may give explicit equations using trigonometric functions). The is no curve which is the *graph of a function* $u = k(x)$ and which stays in the non-shaded area while projecting onto $\mathbb{R} \setminus \{0\}$.

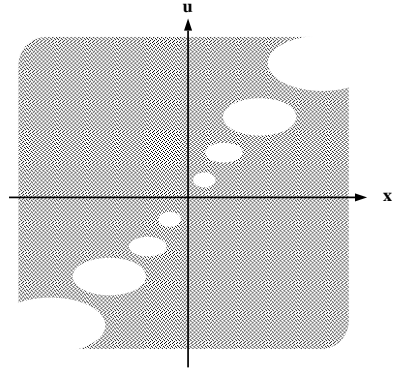

Figure 3: Infinitely many discontinuities required

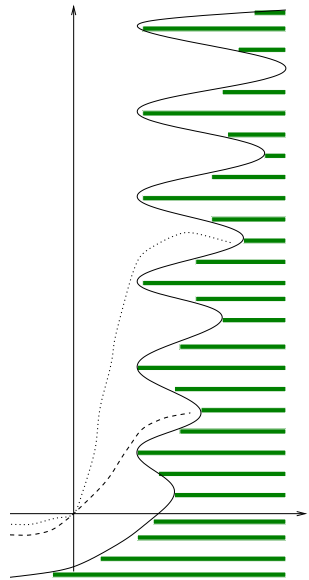

Figure 4: Semiglobal but no global continuous stabilizers

figure shows graphs of two stabilizers, which work with respect to two possible intervals in $\mathbb{R}$, and for any such interval one may clearly find such a (continuous) curve; however, there

**Time-Varying Feedback**

Although our interest here is primarily in time-invariant feedback $u = k(x)$, it is worth pointing out that it is often possible to overcome obstructions to regular feedback stabilization by means of the use of *time-varying* feedback $u = k(t, x)$. A general result in that direction was proved in [39], which established that every one-dimensional system can be stabilized using $k(t, x)$ continuous (it is not difficult to see that the proof can be adapted to obtain $k(t, x)$ periodic in $t$, for each $x$). The intuitive idea of the construction can be perhaps understood by looking at Figure 2.1. Two feedback laws, regular but not stabilizing (because their graphs intersect the shaded area which indicates the complement of $\mathcal{O}$) are shown in the figure. The "bad" set $B = B_t := \{x \mid xf(x, k(t,x)) > 0\}$ where the closed-loop behavior is undesirable moves with time $t$. The idea is to oscillate now between laws whose graphs are between these two graphs, moving slowly in the "$\leftarrow$" direction and faster to the left. The effect will be that states to the right of $B$ move to the left until they hit $B$ (and if $B$ moves slowly enough, never actually hit while $B$ is moving to the left), but when $B$ moves back to the right, the effect on any given state is very short (and hence only a small backwards movement happens). More recently, a very different construction of time-varying feedbacks was accomplished by Coron, cf. [8], who established a general result valid in any dimension, but restricted to systems with no drift, and by Coron and Rosier, cf. [9], for systems for which smooth clf's (see later) exist. We will not look into time-varying feedback in this paper, however.

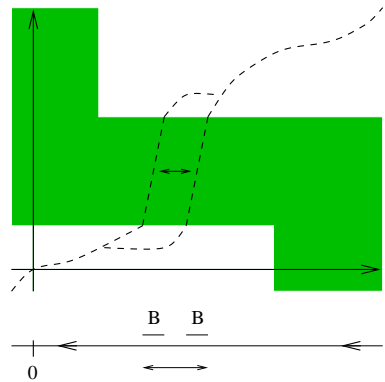

Figure 5: Time-varying feedback

## 2.2 Obstructions

When feedback laws are required to be continuous at the origin, new obstructions arise. The case of systems with $n = m = 1$ is also a good way to introduce this subject. The first observation is that stabilization about the origin (even if just local) means that we must

have, near zero:

$$f(x,k(x)) \begin{cases} > 0 & \text{if } x < 0 \\ < 0 & \text{if } x > 0 \\ = 0 & \text{if } x = 0. \end{cases}$$

In fact, all that we need is that $f(x_1, k(x_1)) < 0$ for some $x_1 > 0$ and $f(x_2, k(x_2)) > 0$ for some $x_2 < 0$. This guarantees, via the intermediate-value theorem that, if $k$ is continuous, the projection

$$(-\varepsilon, \varepsilon) \to \mathbb{R}, \quad x \mapsto f(x, k(x))$$

is onto a neighborhood of zero, for each $\varepsilon > 0$, see Figure 6. It follows, in particular, that

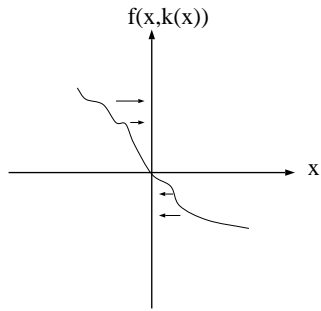

Figure 6: Onto projection, for case $n = m = 1$

the image of

$$(-\varepsilon, \varepsilon) \times (-\varepsilon, \varepsilon) \to \mathbb{R}, \quad (x, u) \mapsto f(x, u)$$

also contains a neighborhood of zero, for any $\varepsilon > 0$ (that is, the map $(x, u) \mapsto f(x, u)$ is open at zero). This last property is intrinsic, being stated in terms of the original data $f(x, u)$ and not depending upon the feedback $k$. Brockett's condition, to be described next, is a far-reaching generalization of this argument; in its proof, degree theory replaces the use of the intermediate value theorem.

### Logical Decisions are Often Necessary

If there are *global obstacles* in the state-space (that is, if the state-space is a proper subset of $\mathbb{R}^n$), discontinuities in feedback laws cannot in general be avoided. Even if it is in principle possible to reach the origin, it may not be possible to find a regular feedback stabilizer. Actually, this is fairly obvious, and is illustrated in intuitive terms by Figure 7. We think of the position of the (immobile) cat as the origin. In deciding in which way to move, as a function of its current position, the dog must at some point in the state-space make a discontinuous decision: move to the left or to the right of the obstacle (represented by the shaded rectangle)? Formally, this setup can be modeled as a problem in which the state-space is the complement in $\mathbb{R}^2$ of the obstacle, and the fact that discontinuities are necessary is a particular case of a general fact (a theorem of Milnor's), namely that the domain of attraction of an asymptotically stable vector field must be diffeomorphic to Euclidean space (which the complement of the rectangle is not); see [36] for more on the subject.

The interesting point is that even if, as in this exposition, we assume that states evolve in Euclidean spaces, similar obstructions may arise. These are due not to the topology of the state space, but to "virtual obstacles" implicit in the form of the system equations. These obstacles occur when it is impossible to move *instantaneously* in certain directions, even if it is possible to move *eventually* in every direction ("nonholonomy").

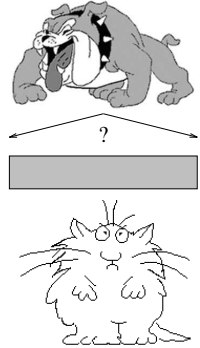

Figure 7: At some point, a discontinuous decision is necessary

**Nonholonomy and Brockett's Theorem**

As an illustration, let us consider a model for the "shopping cart" shown in Figure 8 ("knife-edge" or "unicycle" are other names for this example). The state is given by the orientation $\theta$, together with the coordinates $x_1, x_2$ of the midpoint between the back wheels. The front wheel is a castor, free to rotate. There is a non-slipping constraint on movement: the velocity $(\dot{x}_1, \dot{x}_2)'$ must be parallel to the vector $(\cos\theta, \sin\theta)'$. This leads to the following equations:

$$\begin{aligned} \dot{x}_1 &= u_1 \cos\theta \\ \dot{x}_2 &= u_1 \sin\theta \\ \dot{\theta} &= u_2 \end{aligned}$$

where we may view $u_1$ as a "drive" command and $u_2$ as a steering control; in practice, one would implement these controls by means of differential forces on the two back corners of the cart. The feedback transformation $z_1 := \theta$, $z_2 := x_1 \cos\theta + x_2 \sin\theta$, $z_3 := x_1 \sin\theta - x_2 \cos\theta$, $v_1 := u_2$, and $v_2 := u_1 - u_2 z_3$ brings the system into the system with equations $\dot{z}_1 = v_1$, $\dot{z}_2 = v_2$, $\dot{z}_3 = z_1 v_2$ known as "Brockett's example" or "nonholonomic integrator" (yet another change can bring the third equation into the form $\dot{z}_3 = z_1 v_2 - z_2 v_1$). We view the system as having state space $\mathbb{R}^3$. Although a physically more accurate state space would be the manifold $\mathbb{R}^2 \times \mathbb{S}^1$, the necessary condition to be given is of a local nature, so the global structure is unimportant.

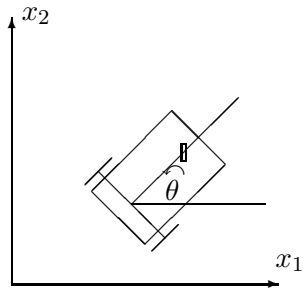

Figure 8: Shopping cart

This system is (obviously) completely controllable (in any case, controllability can be checked using the Lie algebra rank condition, as in e.g. [36], Exercise 4.3.16), and in particular is AC. But we may expect that discontinuities are unavoidable due to the non-slip constraint, which does not allow moving from, for example the position $x_1 = 0$, $\theta = 0$, $x_2 = 1$ in a straight line towards the origin. This is illustrated intuitively by Figure 9 (we use a line to indicate the orientation of the cart). Depending on the initial angle, we may start by rotating

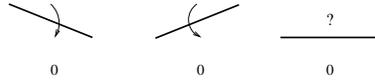

Figure 9: Which way to turn?

either clockwise or counterclockwise (and then moving in a straight line towards 0). But a discontinuous decision is necessary somewhere. This impossibility is established rigorously by means of:

**Theorem A** (Brockett [4]) *If there is a stabilizing feedback which is regular and continuous at zero, then the map $(x, u) \mapsto f(x, u)$ is open at zero.*

The test fails here, since no points of the form $(0, \varepsilon, *)$ belong to the image of the map

$$\mathbb{R}^5 \to \mathbb{R}^3 : (x_1, x_2, \theta, u_1, u_2)' \mapsto f(x, u) = (u_1 \cos \theta, u_1 \sin \theta, u_2)'$$

for $\theta \in (-\pi/2, \pi/2)$, unless $\varepsilon = 0$.

More generally, it is impossible to continuously stabilize any system without drift

$$\dot{x} = u_1 g_1(x) + \ldots + u_m g_m(x) = G(x) u$$

if $m < n$ and $\mathrm{rank}[g_1(0), \ldots, g_m(0)] = m$ (this includes all totally nonholonomic mechanical systems). Indeed, under these conditions, the map $(x, u) \mapsto G(x)u$ cannot contain a neighborhood of zero in its image, when restricted to a small enough neighborhood of zero. Indeed, let us first rearrange the rows of $G$:

$$G(x) \rightsquigarrow \begin{pmatrix} G_1(x) \\ G_2(x) \end{pmatrix}$$

so that $G_1(x)$ is of size $m \times m$ and is nonsingular for all states $x$ that belong to some neighborhood $N$ of the origin. Then,

$$\begin{pmatrix} 0 \\ a \end{pmatrix} \in \mathrm{Im}\, [N \times \mathbb{R}^m \to \mathbb{R}^n : (x, u) \mapsto G(x) u] \;\Rightarrow\; a = 0$$

(since $G_1(x) u = 0 \Rightarrow u = 0 \Rightarrow G_2(x) u = 0$ too).

If the condition $\mathrm{rank}[g_1(0), \ldots, g_m(0)] = m$ is violated, we cannot conclude a negative result. For instance, the system $\dot{x}_1 = x_1 u$, $\dot{x}_2 = x_2 u$ has $m = 1 < 2 = n$ but it can be stabilized by means of the feedback law $u = -(x_1^2 + x_2^2)$.

Observe that for linear systems, Brockett's condition says that

$$\operatorname{rank}[A, B] = n$$

which is the Hautus controllability condition (see e.g. [36], Lemma 3.3.7) at the zero mode.

**Idea of the Proof**

One may prove Brockett's condition in several ways. A proof based on degree theory is probably easiest, and proceeds as follows (for details see for instance [36], Section 5.9). The basic fact, due to Krasnosel'ski, is that if the system $\dot{x} = F(x) = f(x, k(x))$ has the origin as an asymptotically stable point and $F$ is regular (since $k$ is), then the degree (index) of $F$ with respect to zero is $(-1)^n$, where $n$ is the system dimension. In particular, the degree is also nonzero with respect to points $p$ near enough 0, which means that the equation $F(x) = p$ can be solved for small $p$, and hence $f(x, u) = p$ can be solved as well. The proof that the degree is $(-1)^n$ follows by exhibiting a homotopy, namely

$$F_t(x^0) = \frac{1}{t}\left[ x\left(\frac{t}{1-t}, x^0\right) - x^0 \right],$$

between $F_0 = F$ and $F_1(x) = -x$, and noting that the degree of the latter is obviously $(-1)^n$. An alternative (and Brockett's original) proof uses Lyapunov functions. Asymptotic stability implies the existence of a smooth Lyapunov function $V$ for $\dot{x} = F(x) = f(x, k(x))$, so, on the boundary $\partial B$ of a sublevel set $B = \{x \mid V(x) \leq c\}$ we have that $F$ points towards the interior of $B$, see Figure 10. So for $p$ small, $F(x) - p$ still points to the interior, which means that $B$ is invariant with respect to the perturbed vector field $\dot{x} = F(x) - p$. *Provided that a fixed-point theorem applies to continuous maps* $B \to B$, this implies that $F(x) - p$ must vanish somewhere in $B$, that is, the equation $F(x) = p$ can be solved. (Because, for each small $h > 0$, the time-$h$ flow $\phi$ of $F - p$ has a fixed point $x_h \in B$, i.e. $\phi(h, x_h) = x_h$, so picking a convergent subsequence $x_h \to \bar{x}$ gives that $0 = \frac{\phi(h, x_h) - x_h}{h} \to F(\bar{x}) - p$.) A fixed point theorem can indeed be applied, because $B$ is a retract of $\mathbb{R}^n$ (use the flow itself); note that this argument gives a weaker conclusion than the degree condition.

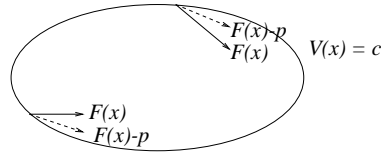

Figure 10: Perturbations of $F$ still point inside $B$

**Another Example**

There is another nice example of these ideas, which will also be useful later when illustrating Lyapunov techniques (see [2]). It is closely related to the shopping-cart example; in fact, it arises when we control the cart by this procedure: first, we rotate the cart until tangent to, either a circle centered on the $x_2$-axis and tangent to the $x_1$-axis (see Figure 11), or the $x_1$-axis if we started there; next, we move only with velocities tangential to this circle (steering so as to maintain invariance of the circle). In summary, one obtains a system with state-space $\mathbb{R}^2$, input space $\mathbb{R}$, and equations

$$\dot{x} = g(x)u, \text{ where } g(x) = \begin{pmatrix} x_1^2 - x_2^2 \\ 2x_1x_2 \end{pmatrix}.$$

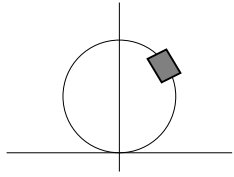

Figure 11: Cart tangent to circle

The vector field $g$ and typical orbits of $g$ are shown in Figure 12. In this system, all motions along the integral curves of $g$ are allowed (clockwise, counterclockwise, or staying at one point). The Brockett condition is satisfied; this is clear if one views the system as a system with complex state space (write $z = x_1 + ix_2$) with equations:

$$\dot{z} = z^2 u.$$

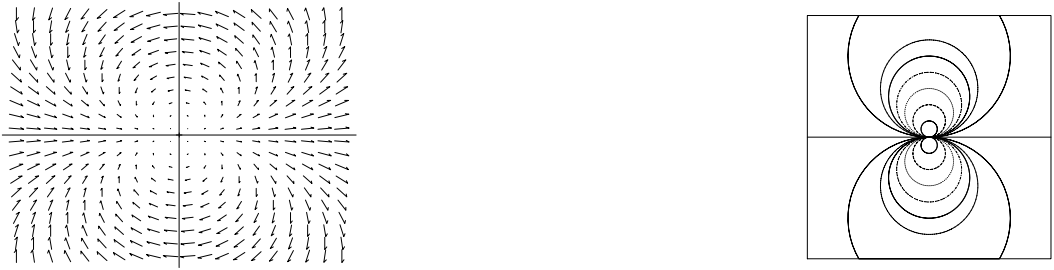

Figure 12: $(x_1^2 - x_2^2)\frac{\partial}{\partial x_1} + 2x_1x_2\frac{\partial}{\partial x_2}$ with typical integral manifolds

The map $(x, u) \mapsto g(x)u$ is $(z, u) \mapsto z^2 u$, which is clearly onto any neighborhood of zero, even when restricted to any neighborhood of $z=u=0$ (take square roots). However, the degree condition fails, as for any continuous feedback one would have $k(x) \neq 0$ for all $x \neq 0$ (otherwise, an equilibrium results), so one obtains degree 2 for $f(x, k(x)) = g(x)k(x)$, and not $(-1)^2 = 1$.

Actually, for this example it is easy to see directly that no regular stabilizing feedback can exist. One proof is by noticing that circles would have to be invariant, so motions restricted to them would result in globally asymptotically stable vector fields in manifolds not diffeomorphic to Euclidean space, contradicting Milnor's theorem. Another proof is even easier: take any path $\gamma : [0, 1] \to \mathbb{R}^2$ with $\gamma(0) = (1, 0)'$ and $\gamma(1) = (-1, 0)'$ and so that $\gamma(r) \neq 0$ for all $r$. Consider the function $\kappa(r) := k(\gamma(r))$, which is continuous if $k$ was assumed continuous. Since $\kappa(0) < 0$ (because, otherwise, trajectories of $\dot{x} = g(x)k(x)$ starting at $(1, 0)'$ cannot converge to zero, since the positive axis is invariant and $g$ points towards the right there) and $\kappa(1) > 0$ (analogous argument), it follows that $\kappa(r_0) = 0$ for some $r_0 \in (0,1)$. Therefore, the point $\bar{x} = \gamma(r_0) \neq 0$ is such that $k(\bar{x}) = 0$ and is therefore an equilibrium point, contradicting the fact that $u = k(x)$ is a stabilizer. We shall return to this example when discussing control-Lyapunov functions.

## 2.3 Control-Lyapunov Functions

The method of control-Lyapunov functions ("clf's") provides a powerful tool for studying stabilization problems, both as a basis of theoretical developments and as a method for actual feedback design.

Before discussing clf's, let us quickly review the classical concept of Lyapunov functions, through a simple example. Consider first a damped spring-mass system $\ddot{y} + \dot{y} + y = 0$, or, in state-space form with $x_1 = y$ and $x_2 = \dot{y}$, $\dot{x}_1 = x_2$, $\dot{x}_2 = -x_1 - x_2$. One way to verify global asymptotic stability of the equilibrium $x = 0$ is to pick the (Lyapunov) function $V(x_1, x_2) := \frac{3}{2} x_1^2 + x_1 x_2 + x_2^2$, and observe that $\nabla V(x).f(x) = -|x|^2 < 0$ if $x \neq 0$, which means that

$$\frac{dV(x(t))}{dt} = -|x(t)|^2 < 0$$

along all nonzero solutions, and thus the energy-like function $V$ decreases along all trajectories, which, since $V$ is a nondegenerate quadratic form, implies that $x(t)$ decreases, and in fact $x(t) \to 0$. Of course, in this case one could compute solutions explicitly, or simply note that the characteristic equation has all roots with negative real part, but Lyapunov functions are a general technique. (In fact, the classical converse theorems of Massera and Kurzweil show that, whenever a system is GAS, there always exists a smooth Lyapunov function $V$.)

Now let us modify this example to deal with a control system, and consider a forced (but undamped) harmonic oscillator $\ddot{x} + x = u$, i.e. $\dot{x}_1 = x_2$, $\dot{x}_2 = -x_1 + u$. The damping feedback $u = -x_2$ stabilizes the system, but let us pretend that we do not know that. If we take the same $V$ as before, now the derivatives along trajectories are, using "$\dot{V}(x,u)$" to denote $\nabla V(x).f(x,u)$ and omitting arguments $t$ in $x(t)$ and $u(t)$:

$$\dot{V}(x,u) = -x_1^2 + x_1 x_2 + x_2^2 - (x_1 + 2x_2)u\,.$$

This expression is affine in $u$. Thus, if $x$ is a state such that $x_1 + 2x_2 \neq 0$, then we may pick a control value $u$ (which depends on this current state $x$) such that $\dot{V} < 0$. On the other hand, if $x_1 + 2x_2 = 0$, then the expression reduces to $\dot{V} = -5x_2^2$ (for any $u$), which is negative unless $x_2$ (and hence also $x_1 = -2x_2$) vanishes.

In conclusion, for each $x \neq 0$ there is some $u$ so that $\dot{V}(x,u) < 0$. This is, except for some technicalities to be discussed, the characterizing property of control-Lyapunov functions. For any given compact subset $B$ in $\mathbb{R}^n$, we now pick some compact subset $\mathcal{U}_0 \subset \mathcal{U}$ so that

$$\forall x \in B, x \neq 0, \quad \exists u \in \mathcal{U}_0 \quad \text{such that} \quad \dot{V}(x,u) < 0\,. \tag{5}$$

In principle, then, we could then stabilize the system, for states in $B$, by using the steepest descent feedback law:

$$k(x) := \operatorname*{argmin}_{u \in \mathcal{U}_0} \nabla V(x) \cdot f(x,u) \tag{6}$$

("argmin" means "pick any $u$ at which the min is attained"; we restricted $\mathcal{U}$ to be assured that $\dot{V}(x,u)$ attains a minimum). Note that the stabilization problem becomes, in these terms, a set of static nonlinear programming problems: minimize a function of $u$, for each $x$.

Global stabilization is also possible, by appropriately picking $\mathcal{U}_0$ as a function of the norm of $x$; later we discuss a precise formulation.

Control-Lyapunov functions, if understood non-technically as the basic paradigm "look for a function $V(x)$ with the properties that $V(x) \approx 0$ if and only if $x \approx 0$, and so that for each $x \neq 0$ it is possible to decrease $V(x)$ by some control action," constitute a very general approach to control (sometimes expressed in a dual fashion, as maximization of some measure of success). They appear in such disparate areas as A.I. game-playing programs (position evaluations), energy arguments for dissipative systems, program termination (Floyd/Dijkstra "variant"), and learning control ("critics" implemented by neural-networks). More relevantly to this paper, the idea underlies much of modern feedback control design, as illustrated for instance by the books [11, 17, 24, 23, 36].

### Differentiable clf's: Precise Definition

We say that a continuous function

$$V : \mathbb{R}^n \to \mathbb{R}_{\geq 0}$$

is *positive definite* if $V(x) = 0$ only if $x = 0$, and it is *proper* (or "weakly coercive") if for each $a \geq 0$ the set $\{x \mid V(x) \leq a\}$ is compact, or, equivalently, $V(x) \to \infty$ as $|x| \to \infty$ (radial unboundedness). A property which is equivalent to properness and positive definiteness together is:

$$(\exists \, \underline{\alpha}, \overline{\alpha} \in \mathcal{K}_\infty) \quad \underline{\alpha}(|x|) \leq V(x) \leq \overline{\alpha}(|x|) \quad \forall \, x \in \mathbb{R}^n \, . \tag{7}$$

A *differentiable control-Lyapunov function* (clf) is a differentiable function $V : \mathbb{R}^n \to \mathbb{R}_{\geq 0}$ which is proper, positive definite, and *infinitesimally decreasing*, meaning that there exists a positive definite continuous function $W : \mathbb{R}^n \to \mathbb{R}_{\geq 0}$, and there is some $\sigma \in \mathcal{N}$, so that

$$\sup_{x \in \mathbb{R}^n} \min_{|u| \leq \sigma(|x|)} \nabla V(x) \cdot f(x, u) + W(x) \leq 0 \, . \tag{8}$$

This is basically the same as condition (5), with $\mathcal{U}_0 = $ the ball of radius $\sigma(|x|)$ picked as a function of $x$. The main difference is that, instead of saying "$\nabla V(x) \cdot f(x, u) < 0$ for $x \neq 0$" we write $\nabla V(x) \cdot f(x, u) \leq -W(x)$, where $W$ is negative when $x \neq 0$. The two definitions are equivalent, but the "Hamiltonian" version used here is the correct one for the generalizations to be given, to nonsmooth $V$.

**Theorem B** (Artstein [2]) *A control-affine system $\dot{x} = g_0(x) + \sum u_i g_i(x)$ admits a differentiable clf if and only if it admits a regular stabilizing feedback.*

The proof of sufficiency is easy: if there is such a $k$, then the converse Lyapunov theorem, applied to the closed-loop system $F(x) = f(x, k(x))$, provides a smooth $V$ such that

$$L_F V(x) = \nabla V(x) F(x) < 0 \quad \forall x \neq 0 \, .$$

This gives that for all nonzero $x$ there is some $u$ (bounded on bounded sets, because $k$ is locally bounded by definition of feedback) so that $\dot{V}(x, u) < 0$; and one can put this in the form (8).

The necessity is more interesting. The original proof in [2] proceeds by a nonconstructive argument involving partitions of unity, but it is also possible to exhibit explicitly a feedback, written as a function
$$k\left(\nabla V(x) \cdot g_0(x),\, \ldots\,,\, \nabla V(x) \cdot g_m(x)\right)$$
of the directional derivatives of $V$ along the vector fields defining the system (*universal formulas* for stabilization). Taking for simplicity $m = 1$, one such formula is:

$$k(x) \;:=\; -\,\frac{a(x) + \sqrt{a(x)^2 + b(x)^4}}{b(x)} \quad (0 \text{ if } b = 0)$$

where $a(x) := \nabla V(x) \cdot g_0(x)$ and $b(x) := \nabla V(x) \cdot g_1(x)$. (The expression for $k$ is analytic in $a, b$ when $x \neq 0$, because the clf property means that $a(x) < 0$ whenever $b(x) = 0$, see [36] for details.)

Thus, the question of existence of regular feedback, for control-affine systems, reduces to the search for differentiable clf's, and this gives rise to a vast literature dealing with the construction of such $V$'s, see [11, 24, 23, 36] and references therein. Many other theoretical issues are also answered by Artstein's theorem. For example, via Kurzweil's converse theorem one has that the existence of $k$ merely continuous on $\mathbb{R}^n \setminus \{0\}$ suffices for the existence of smooth (infinitely differentiable) $V$, and from here one may in turn find a $k$ which is smooth on $\mathbb{R}^n \setminus \{0\}$. In addition, one may easily characterize the existence of $k$ continuous at zero as well as regular: this is equivalent to the *small control property*: for each $\varepsilon > 0$ there is some $\delta > 0$ so that $0 < |x| < \delta$ implies that $\min_{|u| \leq \varepsilon} \nabla V(x) \cdot f(x, u) < 0$ (if this property holds, the universal formula automatically provides such a $k$). We should note that Artstein provided a result valid for general, not necessarily control-affine systems $\dot{x} = f(x, u)$; however, the obtained "feedback" has values in sets of relaxed controls, and is not a feedback law in the classical sense. Later, we discuss a different generalization.

Differentiable clf's will in general not exist, because of obstructions to regular feedback stabilization. This leads us naturally into the twin subjects of discontinuous feedbacks and non-differentiable clf's.

## 3 Discontinuous Feedback

The previous results and examples show that, in order to develop a satisfactory general theory of stabilization, one in which one proves the implication "asymptotic controllability implies feedback stabilizability," we must allow discontinuous feedback laws $u = k(x)$. But then, a major technical difficulty arises: solutions of the initial-value problem $\dot{x} = f(x, k(x))$, $x(0) = x^0$, interpreted in the classical sense of differentiable functions or even as (absolutely) continuous solutions of the integral equation $x(t) = x^0 + \int_0^t f(x(s), k(x(s)))\,ds$, do not exist in general. The only general theorems apply to systems $\dot{x} = F(x)$ with continuous $F$. For example, there is no solution to $\dot{x} = -\text{sign}\,x$, $x(0) = 0$, where $\text{sign}\,x = -1$ for $x < 0$ and $\text{sign}\,x = 1$ for $x \geq 0$. So one cannot even pose the stabilization problem in a mathematically consistent sense.

There is, of course, an extensive literature addressing the question of discontinuous feedback laws for control systems and, more generally, differential equations with discontinuous right-hand sides. One of the best-known candidates for the concept of solution of (3) is that of a *Filippov solution* [12, 13], which is defined as the solution of a certain differential inclusion with a multivalued right-hand side which is built from $f(x, k(x))$. Unfortunately, there is no hope of obtaining the implication "asymptotic controllability implies feedback stabilizability" if one interprets solutions of (3) as Filippov solutions. This is a consequence of results in [32, 9], which established that the existence of a discontinuous stabilizing feedback in the Filippov sense implies the Brockett necessary conditions, and, moreover, for systems affine in controls it also implies the existence of regular feedback (which we know is in general impossible).

A different concept of solution originates with the theory of discontinuous positional control developed by Krasovskii and Subbotin in the context of differential games in [21], and it is the basis of the new approach to discontinuous stabilization proposed in [7], to which we now turn.

### 3.1 Limits of High-Frequency Sampling

By a *sampling schedule* or *partition* $\pi = \{t_i\}_{i \geq 0}$ of $[0, +\infty)$ we mean an infinite sequence

$$0 = t_0 < t_1 < t_2 < \ldots$$

with $\lim_{i \to \infty} t_i = \infty$. We call

$$\mathbf{d}(\pi) := \sup_{i \geq 0}(t_{i+1} - t_i)$$

the *diameter* of $\pi$. Suppose that $k$ is a given feedback law for system (1). For each $\pi$, the *$\pi$-trajectory starting from $x^0$* of system (3) is defined recursively on the intervals $[t_i, t_{i+1})$, $i = 0, 1, \ldots$, as follows. On each interval $[t_i, t_{i+1})$, the initial state is measured, the control value $u_i = k(x(t_i))$ is computed, and the constant control $u \equiv u_i$ is applied until time $t_{i+1}$; the process is then iterated. That is, we start with $x(t_0) = x^0$ and solve recursively

$$\dot{x}(t) \;=\; f(x(t), k(x(t_i))), \; t \in [t_i, t_{i+1}), \quad i = 0, 1, 2, \ldots$$

using as initial value $x(t_i)$ the endpoint of the solution on the preceding interval. The ensuing $\pi$-trajectory, which we denote as $x_\pi(\cdot, x^0)$, is defined on some maximal nontrivial interval; it may fail to exist on the entire interval $[0, +\infty)$ due to a blow-up on one of the subintervals $[t_i, t_{i+1})$. We say that it is *well defined* if $x_\pi(t, x^0)$ is defined on all of $[0, +\infty)$.

**Definition.** The feedback $k : \mathbb{R}^n \to \mathcal{U}$ *stabilizes* the system (1) if there exists a function $\beta \in \mathcal{KL}$ so that the following property holds: For each

$$0 < \varepsilon < K$$

there exists a $\delta = \delta(\varepsilon, K) > 0$ such that, for every sampling schedule $\pi$ with $\mathbf{d}(\pi) < \delta$, and for each initial state $x^0$ with $|x^0| \leq K$, the corresponding $\pi$-trajectory of (3) is well-defined and satisfies

$$|x_\pi(t, x^0)| \;\leq\; \max\left\{\beta\left(K, t\right), \varepsilon\right\} \quad \forall\, t \geq 0 \,. \tag{9}$$

In particular, we have

$$|x_\pi(t, x^0)| \;\leq\; \max\left\{\beta\left(|x^0|, t\right), \varepsilon\right\} \quad \forall\, t \geq 0 \tag{10}$$

whenever $0 < \varepsilon < |x^0|$ and $\mathbf{d}(\pi) < \delta(\varepsilon, |x^0|)$ (just take $K := |x^0|$).

Observe that the role of $\delta$ is to specify a lower bound on intersampling times. Roughly, one is requiring that
$$t_{i+1} \leq t_i + \theta\left(|x(t_i)|\right)$$
for each $i$, where $\theta$ is an appropriate positive function.

Our definition of stabilization is physically meaningful, and is very natural in the context of sampled-data (computer control) systems. It says in essence that a feedback $k$ stabilizes the system if it drives all states asymptotically to the origin and with small overshoot when using *any fast enough sampling schedule*. A high enough sampling frequency is generally required when close to the origin, in order to guarantee small displacements, and also at infinity, so as to preclude large excursions or even blow-ups in finite time. This is the reason for making $\delta$ depend on $\varepsilon$ and $K$.

This concept of stabilization can be reinterpreted in various ways. One is as follows. Pick any initial state $x^0$, and consider any sequence of sampling schedules $\pi_\ell$ whose diameters $\mathbf{d}(\pi_\ell)$ converge to zero as $\ell \to \infty$ (for instance, constant sampling rates with $t_i = i/\ell$, $i = 0, 1, 2, \ldots$). Note that the functions $x_\ell := x_{\pi_\ell}(\cdot, x^0)$ remain in a bounded set, namely the ball of radius $\beta(|x^0|, 0)$ (at least for $\ell$ large enough, for instance, any $\ell$ so that $\mathbf{d}(\pi_\ell) < \delta(|x^0|/2, |x^0|)$). Because $f(x, k(x))$ is bounded on this ball, these functions are equicontinuous, and (Arzela-Ascoli's Theorem) we may take a subsequence, which we denote again as $\{x_\ell\}$, so that $x_\ell \to x$ as $\ell \to \infty$ (uniformly on compact time intervals) for some absolutely continuous (even Lipschitz) function $x : [0, \infty) \to \mathbb{R}^n$. *We may think of any limit function $x(\cdot)$ that arises in this fashion as a generalized solution of the closed-loop equation (3).* That is, generalized solutions are the limits of trajectories arising from arbitrarily high-frequency sampling when using the feedback law $u = k(x)$. Generalized solutions, for a given initial state $x^0$, may not be unique – just as may happen with continuous but non-Lipschitz feedback – but there is always existence, and, moreover, for any generalized solution, $|x(t)| \leq \beta(|x^0|, t)$ for all $t \geq 0$. This is precisely the defining estimate for the GAS property. Moreover, if $k$ happens to be regular, then the unique solution of $\dot{x} = f(x, k(x))$ in the classical sense is also the unique generalized solution, so we have a reasonable extension of the concept of solution. (This type of interpretation is somewhat analogous, at least in spirit, to the way in which "relaxed" controls are interpreted in optimal trajectory calculations, namely through high-frequency switching of approximating regular controls.)

**Remark.** The definition of stabilization was given in [7] in a slightly different form. It was required there that there exist for each $R > 0$ a number $M(R) > 0$, with $\lim_{R \searrow 0} M(R) = 0$, and, for each $0 < r < R$, numbers $\delta_0(r, R) > 0$ and $T(r, R) \geq 0$, such that the following property holds: for each sampling schedule $\pi$ with $\mathbf{d}(\pi) < \delta_0(r, R)$, each $x^0$ with $|x^0| \leq R$, and each $t \geq T(r, R)$, it holds that $|x_\pi(t, x^0)| \leq r$, and, in addition, $|x_\pi(t, x^0)| \leq M(R)$ for all $t \geq 0$. This definition is equivalent to the definition that we gave, based on an estimate of the type (9). See Section 5.1 for a proof.

## 3.2 Stabilizing Feedbacks Exist

In the paper [7], the following result was proven by Clarke, Ledyaev, Subbotin, and the author:

**Theorem C** *The system (1) admits a stabilizing feedback if and only if it is asymptotically controllable.*

Necessity is clear. The sufficiency statement is proved by construction of $k$, and is based on the following ingredients:

- Existence of a nonsmooth control-Lyapunov function $V$.
- Regularization on shells of $V$.
- Pointwise minimization of a Hamiltonian for the regularized $V$

In order to sketch this construction, we start by quickly reviewing a basic concept from nonsmooth analysis.

**Proximal Subgradients**

Let $V$ be any continuous function $\mathbb{R}^n \to \mathbb{R}$ (or even, just lower semicontinuous and with extended real values). A *proximal subgradient* of $V$ at the point $x \in \mathbb{R}^n$ is any vector $\zeta \in \mathbb{R}^n$ such that, for some $\sigma > 0$ and some neighborhood $\mathcal{O}$ of $x$,

$$V(y) \geq V(x) + \zeta \cdot (y - x) - \sigma^2 |y - x|^2 \quad \forall y \in \mathcal{O}.$$

In other words, proximal subgradients are the possible gradients of supporting quadratics at the point $x$. The set of all proximal subgradients at $x$ is denoted $\partial_P V(x)$. For example (see Figure 13), the function $V(x) = |x|$ admits any $\zeta \in [-1, 1]$ as a proximal subgradient at $x = 0$ (elsewhere, $\partial_P V(x) = \{\nabla V(x)\}$), while for $V(x) = -|x|$ we have $\partial_P V(0) = \emptyset$ (because there are no possible quadrics that fit inside the graph and touch the corner).

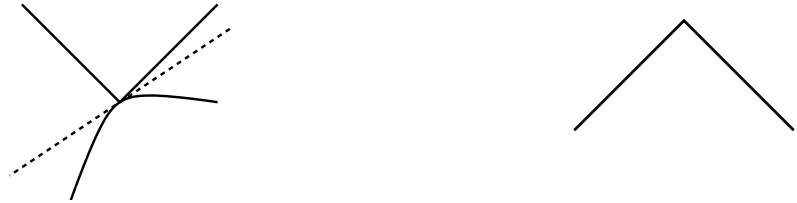

Figure 13: $\partial_P |x|(0) = [-1, 1]$     $\partial_P(-|x|)(0) = \emptyset$

**Nonsmooth Control-Lyapunov Functions**

A *continuous* (but not necessarily differentiable) $V : \mathbb{R}^n \to \mathbb{R}_{\geq 0}$ is a *control-Lyapunov function* (clf) if it is proper, positive definite, and infinitesimally decreasing in the following generalized sense: there exist a positive definite continuous $W : \mathbb{R}^n \to \mathbb{R}_{\geq 0}$ and a $\sigma \in \mathcal{N}$ so that

$$\sup_{x \in \mathbb{R}^n} \max_{\zeta \in \partial_P V(x)} \min_{|u| \leq \sigma(|x|)} \zeta \cdot f(x, u) + W(x) \leq 0. \tag{11}$$

This is the obvious generalization of the differentiable case in (8); we are still asking that one should be able to make $\nabla V(x) \cdot f(x, u) < 0$ by an appropriate choice of $u = u_x$, for each $x \neq 0$, except that now we replace $\nabla V(x)$ by the proximal subgradient set $\partial_P V(x)$. An equivalent property is to ask that $V$ be a viscosity supersolution of the corresponding Hamilton-Jacobi-Bellman equation.

In the paper [34], the following result was proven by the author:

**Theorem D** *The system (1) is asymptotically controllable if and only if it admits a continuous clf.*

Not surprisingly, the proof is based on first constructing an appropriate $W$, and then letting $V$ be the optimal cost (Bellman function) for the problem $\min \int_0^\infty W(x(s))\,ds$. However, some care has to be taken to insure that $V$ is continuous, and the cost has to be adjusted in order to deal with possibly unbounded minimizers. Actually, to be precise, the result as stated here is really a restatement (cf. [40], [7]) of the main theorem given in [34]. See Section 5.2 for the details of this reduction.

**Regularization**

Once $V$ is known to exist, the next step in the construction of a stabilizing feedback is to obtain Lipschitz approximations of $V$. For this purpose, one considers the Iosida-Moreau inf-convolution of $V$ with a quadratic function:

$$V_\alpha(x) := \inf_{y \in \mathbb{R}^n}\left[V(y) + \frac{1}{2\alpha^2}|y-x|^2\right]$$

where the number $\alpha > 0$ is picked constant on appropriate regions. One has that $V_\alpha(x) \nearrow V(x)$, uniformly on compacts. Since $V_\alpha$ is locally Lipschitz, Rademacher's Theorem insures that $V_\alpha$ is differentiable almost everywhere. The feedback $k$ is then made equal to a pointwise minimizer $k_\alpha$ of the Hamiltonian, at the points of differentiability (compare with (6) for the case of differentiable $V$):

$$k_\alpha(x) := \operatorname*{argmin}_{u \in \mathcal{U}_0} \nabla V_\alpha(x) \cdot f(x,u),$$

where $\alpha$ and the compact $\mathcal{U}_0 = \mathcal{U}_0(\alpha)$ are chosen constant on certain compacts and this choice is made in between level curves, see Figure 14. The critical fact is that $V_\alpha$ is itself a clf for the original system, at least when restricted to the region where it is needed. More precisely, on each shell of the form

$$C = \{x \in \mathbb{R}^n \mid r \leq |x| \leq R\},$$

there are positive numbers $m$ and $\alpha_0$ and a compact subset $\mathcal{U}_0$ such that, for each $0 < \alpha \leq \alpha_0$, each $x \in C$, and every $\zeta \in \partial_P V_\alpha(x)$,

$$\min_{u \in \mathcal{U}_0} \zeta \cdot f(x,u) + m \leq 0.$$

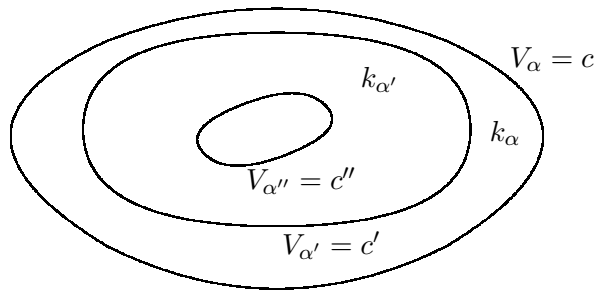

Figure 14: $k = k_\alpha$ on $\{x \mid V_\alpha(x) \leq c, V_{\alpha'}(x) > c'\}$

See Section 5.3.

Actually, this description is oversimplified, and the proof is a bit more delicate. One must define, on appropriate compact sets

$$k(x) := \operatorname*{argmin}_{u \in \mathcal{U}_0} \zeta_\alpha(x) \cdot f(x, u),$$

where $\zeta_\alpha(x)$ is carefully chosen. At points $x$ of nondifferentiability, $\zeta_\alpha(x)$ is not a proximal subgradient of $V_\alpha$, since $\partial_P V_\alpha(x)$ may well be empty. One uses, instead, the fact that $\zeta_\alpha(x)$ happens to be in $\partial_P V(x')$ for some $x' \approx x$.

**An Example**

As a simple example, we consider the system that was obtained from the two-dimensional reduction of the "shopping cart" problem, cf. Figure 12. For this example, continuous stabilization is not possible, and so no differentiable clf's can exist. On the other hand, the system is AC, so one can stabilize it using the techniques just described. A clf for this problem was obtained in [25]:

$$V(x_1, x_2) = \frac{x_1^2 + x_2^2}{\sqrt{x_1^2 + x_2^2} + |x_1|}$$

( 0 if $x_1 = x_2 = 0$) and its level sets are as shown in Figure 15. Note that the nonsmoothness happens exactly on the $x_2$ axis (one has an empty subgradient at those points), and the clf inequality follows from the fact that

$$\inf_u \nabla V(x_1, x_2) \cdot f((x_1, x_2), u) \leq -\frac{1}{2}(x_1^2 + x_2^2)$$

at points with $x_1 \neq 0$ (the proximal subgradient set is empty otherwise). The stabilizing feedback that results is the obvious one: if to the right of the $x_2$ axis, move clockwise, if to the left counterclockwise, and make an arbitrary decision (this arbitrariness corresponds to the choice of "$\zeta_\alpha(x)$" in the theory) on the $x_2$-axis.

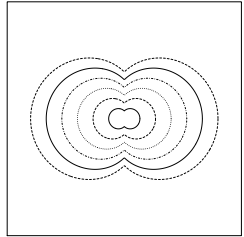

Figure 15: Clf Level Sets for system in Figure 12

## 4   Sensitivity to Small Measurement Errors

We have seen that every asymptotically controllable system admits a feedback stabilizer $k$, generally discontinuous, which renders the closed-loop system $\dot{x} = f(x, k(x))$ GAS. On the other hand, one of the main reasons for using feedback is to deal with uncertainty, and one possible source of uncertainty are measurement errors in state estimation. The use of discontinuous feedback means that undesirable behavior –chattering– may arise. In fact, one of the main reasons for the focus on continuous feedback is precisely in order to avoid such behaviors. Thus, we turn now to an analysis of the effect of measurement errors.

Suppose first that $k$ is a continuous function of $x$. Then, if the error $e$ is small, using the control $u' = k(x + e)$ instead of $u = k(x)$ results in behavior which remains close to the intended one, since $k(x + e) \approx k(x)$; moreover, if $e \ll x$ then stability is preserved. This property of robustness to small errors when $k$ is continuous can be rigorously established by means of a Lyapunov proof, based on the observation that, if $V$ is a Lyapunov function for the closed-loop system, then continuity of $f(x, k(x + e))$ on $e$ means that

$$\nabla V(x) \cdot f(x, k(x+e)) \approx \nabla V(x) \cdot f(x, k(x)) < 0.$$

Unfortunately, when $k$ is not continuous, this argument breaks down. However, it can be modified so as to avoid invoking continuity of $k$. Assuming that $V$ is continuously differentiable, one can argue that

$$\nabla V(x) \cdot f(x, k(x+e)) \approx \nabla V(x+e) \cdot f(x, k(x+e)) < 0$$

(using the Lyapunov property at the point $x + e$ instead of at $x$). This observation leads to a theorem, formulated below, which says that a discontinuous feedback stabilizer, robust with respect to small observation errors, can be found provided that there is a $\mathcal{C}^1$ clf.

In general, as there are no $\mathcal{C}^1$, but only continuous, clf's, one may not be able to find any feedback law that is robust in this sense. We can see this fact intuitively with an example. Let us take once more the two-dimensional problem illustrated in Figure 12, and let us suppose that we are using the following control law: if to the right of, or exactly on, the $x_2$ axis, move clockwise, and if to the left move counterclockwise. See Figure 16, which indicates what happens on any circle. The main point that we wish to make is that *this feedback law*

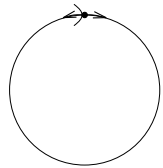

Figure 16: Feedback $k$ on a typical integral manifold

*is extremely sensitive to measurement errors.* Indeed, if the true state $x$ is slightly to the left of the top point, but we mistakenly believe it to be to the right, we use a clockwise motion, in effect bringing the state towards the top, instead of downwards towards the target (the origin); see Figure 17. It is clear that, if we are unlucky enough to consistently make measurement errors that place us on the opposite side of the true state, the result may well be an oscillation (chattering) around the top.

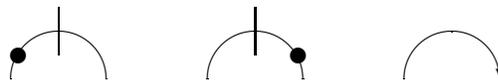

Figure 17: True state; measured state; erroneous motion commanded

We might say, then, that the discontinuous feedback laws which are guaranteed to always exist by the general theorem in [7] lead to *fussy control*\* –not to be confused, of course, with "fuzzy" control.

There are many well-known techniques for avoiding chattering, and a very common one is the introduction of deadzones where no action is taken. Indeed, in the above example, one may adopt the following modified control strategy: stay in the chosen direction (even if it might be "wrong") for some minimal time, until we are guaranteed to be far enough from the discontinuity; only after this minimal amount of time, we sample again. At this point, we know for sure on which side of the top we are. (This assumes that we have an upper bound on the magnitude of the error. Also, of course, observation errors when close to the origin will mean that we can only expect "practical" stability, meaning that we cannot be assured of convergence to the origin, but merely of convergence to a neighborhood of the origin whose size depends on the size of the observation errors.)

Such a control strategy is not a pure "continuous time" one, in that a minimum intersample time is required. It can be interpreted, rather, as constructing a *hybrid* (different time scales needed) and *dynamic* (requiring memory) controller. The paper [26] proves a general result showing the possibility of stabilization of every AC system, using an appropriate definition of general hybrid dynamic controllers. The controller given there incorporates an internal model of system. It compares, at appropriate sampling times, the state predicted by the internal model with the –noisy– observations of the state; whenever these differ substantially, a "reseting" is performed in the state of the controller.

Actually, already the feedback constructed in [7], with no modifications needed, can always be used in a manner robust with respect to small observation errors, using the idea illustrated with the circle of not sampling again for some minimal period. Roughly speaking, the general idea is as follows.

Suppose that the true current state, let us say at time $t = t_i$, is $x$, but that the controller uses $u = k(\tilde{x})$, where $\tilde{x} = x + e$, and $e$ is small. Call $x'$ the state that results at the next sampling time, $t = t_{i+1}$. By continuity of solutions on initial conditions, $|x' - \tilde{x}'|$ is also small, where $\tilde{x}'$ is the state that would have resulted from applying the control $u$ if the true state had been $\tilde{x}$. By continuity, it follows that $V_\alpha(x) \approx V_\alpha(\tilde{x})$ and also $V_\alpha(x') \approx V_\alpha(\tilde{x}')$. On the other hand, the construction in [7] provides that $V_\alpha(\tilde{x}') < V_\alpha(\tilde{x}) - d(t_{i+1} - t_i)$, where $d$ is some positive constant (this is valid while we are far from the origin). Hence, if $e$ is sufficiently small compared to the intersample time $t_{i+1} - t_i$, it will necessarily be the case that $V_\alpha(x')$ must also be smaller than $V_\alpha(x)$. See Figure 18. Thinking of $V_\alpha$ as a Lyapunov function, this means that $x'$ is made "smaller", even if the wrong control $u = k(\tilde{x})$, rather than $k(x)$, is applied.

This discussion may be formalized in several ways. We limit ourselves here to a theorem assuring semiglobal practical stability (i.e., driving all states in a given compact set of initial conditions into a specified neighborhood of zero). For any sampling schedule $\pi$, we denote

$$\underline{\mathbf{d}}(\pi) := \inf_{i \geq 0}(t_{i+1} - t_i).$$

If $e : [0, \infty) \to \mathbb{R}^n$ is any function ($e(t)$ is to be thought of as the state estimation error at time $t$), $k$ is a feedback law, $x^0 \in \mathbb{R}^n$, and $\pi$ is a sampling schedule, we define the solution of

$$\dot{x} \;=\; f(x, k(x+e)), \quad x(0) = x^0 \tag{12}$$

as earlier, namely, recursively solving

$$\dot{x}(t) \;=\; f(x(t), k(x(t_i) + e(t_i)))$$

with initial condition $x(t_i)$ on the intervals $[t_i, t_{i+1}]$. The feedback stabilizer that was constructed in [7] was defined by patching together feedback laws, denoted there as $k_\nu$ (where $\nu = (\alpha, r, R)$ is a triple of positive numbers, with $r < R$). We use these same feedbacks in the statement of the result.

---

\*Fussy (adjective): "...requiring ...close attention to details" (Webster).

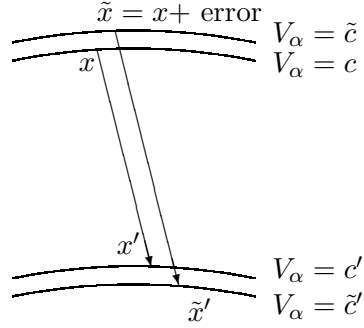

Figure 18: $t_{i+1} \gg t_i \;\Rightarrow\; \tilde{c}' \ll \tilde{c} \;\Rightarrow\; c' < c$

**Theorem E** *Suppose that the system (1) is asymptotically controllable. Then, there exists a function $\Gamma \in \mathcal{K}_\infty$ with the following property. For each $0 < \varepsilon < K$, there is a feedback of the type $k_\nu$, and there exist positive $\delta = \delta(\varepsilon, K)$, $\kappa = \kappa(\varepsilon, K)$, and $T = T(\varepsilon, K)$, such that, for each sampling schedule $\pi$ with $\mathbf{d}(\pi) \leq \delta$, each $e : [0, \infty) \to \mathbb{R}^n$ so that*

$$|e(t)| \leq \kappa \underline{\mathbf{d}}(\pi) \quad \forall t \geq 0\,,$$

*and each $x^0$ with $|x^0| \leq K$, the solution of the noisy system (12) satisfies*

$$|x(t)| \;\leq\; \Gamma(K) \quad \forall t \geq 0$$

*and*

$$|x(t)| \;\leq\; \varepsilon \quad \forall t \geq T\,.$$

See Section 5.4.

### 4.1 A Necessary Condition

Theorem E insures that stabilization is possible if we sample "just right" (not too slow, so as to preserve stability, but also not too fast, so that observation errors do not cause chattering). It leaves open the theoretical question of precisely under what conditions is it possible to find a state feedback law which is robust with respect to small observation errors and which, on the other hand, is a continuous-time feedback, in the sense of arbitrarily fast sampling. The discussion preceding Theorem E suggests that this objective cannot always be met (e.g., for the circle problem), but that the existence of a $\mathcal{C}^1$ clf might be sufficient for guaranteeing that it can. Indeed, this is what happens, as was proved in the recent paper [27]. We next present the main result from that paper.

We consider systems

$$\dot{x}(t) \;=\; f(x(t), k(x(t) + e(t)) + d(t)) \tag{13}$$

in which there are observation errors as well as, now, possible actuator errors $d(\cdot)$. (Actually, robustness to just small actuator errors is not a serious issue; the original paper [7] showed that the feedback laws obtained there stabilize even in the presence of such errors, or even model errors.) We assume that actuator errors $d(\cdot) : [0, \infty) \to \mathcal{U}$ are Lebesgue measurable and locally essentially bounded, and that observation errors $e(\cdot) : [0, \infty) \to \mathbb{R}^n$ are locally bounded.

We wish to define what it means for a feedback $k$ to stabilize (13); roughly, we will ask that, for some $\beta \in \mathcal{KL}$, and for some "tolerance" function $\theta : \mathbb{R}_{\geq 0} \to \mathbb{R}_{\geq 0}$, we have an estimate $|x(t, x^0)| \leq \beta(|x^0|, t)$ provided that $|e(t)| \leq \theta(|x(t)|)$ and $|d(t)| \leq \theta(|x(t)|)$ for all $t$ (small enough errors). However, the definition is somewhat complicated by the need to appropriately choose sampling frequencies.

We define solutions of (13), for each sampling schedule $\pi$, in the usual manner, i.e., solving recursively on the intervals $[t_i, t_{i+1})$, $i = 0, 1, \ldots$, the differential equation

$$\dot{x}(t) \;=\; f(x(t), \, k(x(t_i) + e(t_i)) + d(t)) \tag{14}$$

with $x(0) = x^0$. We write $x(t) = x_\pi(t, x^0, d, e)$ for the solution, and say it is *well-defined* if it is defined for all $t \geq 0$.

**Definition.** The feedback $k : \mathbb{R}^n \to \mathcal{U}$ *stabilizes* the system (13) if there exists a function $\beta \in \mathcal{KL}$ so that the following property holds: For each

$$0 < \varepsilon < K$$

there exist $\delta = \delta(\varepsilon, K) > 0$ and $\eta = \eta(\varepsilon, K)$ such that, for every sampling schedule $\pi$ with $\mathbf{d}(\pi) < \delta$, each initial state $x^0$ with $|x^0| \leq K$, and each $e, d$ such that $|e(t)| \leq \eta$ for all $t \geq 0$ and $|d(t)| \leq \eta$ for almost all $t \geq 0$, the corresponding $\pi$-trajectory of (13) is well-defined and satisfies

$$|x_\pi(t, x^0, d, e)| \;\leq\; \max\{\beta(K, t)\, , \, \varepsilon\} \quad \forall\, t \geq 0\,. \tag{15}$$

In particular, taking $K := |x^0|$, one has that

$$|x_\pi(t, x^0, d, e)| \;\leq\; \max\{\beta(|x^0|, t)\, , \, \varepsilon\} \quad \forall\, t \geq 0$$

whenever $0 < \varepsilon < |x^0|$, $\mathbf{d}(\pi) < \delta(\varepsilon, |x^0|)$, and for all $t$, $|e(t)| \leq \eta(\varepsilon, |x^0|)$, and $|d(t)| \leq \eta(\varepsilon, |x^0|)$.

The main result in [27] is as follows.

**Theorem F** *There is a feedback which stabilizes the system (13) if and only if there is a $\mathcal{C}^1$ clf for the unperturbed system (1).*

This result is somewhat analogous to a result obtained, for classical solutions, by Hermes in [14]; see also [13].

It is interesting to note that, as a corollary of Artstein's Theorem, for control-affine systems $\dot{x} = g_0(x) + \sum u_i g_i(x)$ we may conclude that if there is a discontinuous feedback stabilizer that is robust with respect to small noise, then there is also a regular one, and even one that is smooth on $\mathbb{R}^n \setminus \{0\}$.

For non control-affine systems, however, there may exist a discontinuous feedback stabilizer that is robust with respect to small noise, yet there is no regular feedback. For example, consider the following system with $n = 3$ and $m = 1$:

$$\begin{aligned} \dot{x}_1 &= u_2 u_3 \\ \dot{x}_2 &= u_1 u_3 \\ \dot{x}_3 &= u_1 u_2 \,. \end{aligned}$$

Here, there is a $\mathcal{C}^1$ clf, namely the squared norm $(x_1^2 + x_2^2 + x_3^2)$, but there is no possible regular feedback stabilizer, since Brockett's condition fails because points $(0, \neq 0, \neq 0)$ cannot be in the image of $(x, u) \mapsto f(x, u)$. (Because this is a homogeneous system with no drift, Brockett's condition rules out even feedbacks that are not continuous at the origin, see [27] for a remark to that effect.)

The sufficiency part of Theorem F proceeds by taking a pointwise minimization of the Hamiltonian, for a given $\mathcal{C}^1$ clf, i.e. $k(x)$ is defined as any $u$ with $|u| \leq \sigma(|x|)$ which minimizes $\nabla V(x) \cdot f(x, u)$. The necessity part is based on the following technical fact: if the perturbed system can be stabilized, then the differential inclusion

$$\dot{x} \in F(x) := \bigcap_{\varepsilon > 0} \overline{\operatorname{co}} f(x, k(x + \varepsilon B))$$

(where $B$ denotes the unit ball in $\mathbb{R}^n$) is strongly asymptotically stable. One may then apply the recent converse Lyapunov theorem of [5] for upper semicontinuous compact convex differential inclusions (which generalized to differential inclusions the theorem from [28] characterizing uniform asymptotic stability of systems with disturbances $\dot{x} = f(x, d)$) to deduce the existence of $V$.

We can now summarize exactly which implications hold. We write "robust" to mean stabilization of the system subject to observation and actuator noise:

$$\begin{array}{ccccc} \mathcal{C}^1 V & \Longleftrightarrow & \exists \text{ robust } k & & \\ \Downarrow & & \Downarrow & & \\ \mathcal{C}^0 V & \Longleftrightarrow & \exists k & \Longleftrightarrow & \text{AC} \end{array}$$

**Remark.** The definition of stabilization of (13) was given in [27] in a slightly different form than here. There, it was required that for each $0 < r < R$ there exist $M = M(R) > 0$ with $\lim_{R \searrow 0} M(R) = 0$, $\delta = \delta(r, R) > 0$, $T = T(r, R) > 0$, and $\eta = \eta(r, R)$, such that, for every partition $\pi$ with $\mathbf{d}(\pi) < \delta$, each initial state with $|x^0| \leq R$, and each $e, d$ such that $|e(t)| \leq \eta$ for all $t \geq 0$ and $|d(t)| \leq \eta$ for almost all $t \geq 0$, the $\pi$-trajectory of

$$\dot{x} = f(x, k(x + e) + d)$$

is defined for all $t \geq 0$ and

$$|x(t)| \leq r \; \forall t \geq T \quad \text{and} \quad |x(t)| \leq M(R) \; \forall t \geq 0 \,.$$

This definition is equivalent to the definition that we just gave. The same proof as in Section 5.1 applies.

## 5 Some Proofs

We fill-in here the proofs of several technical points regarding stabilization.

### 5.1 Estimates for Stabilization

We prove here that if a feedback $k$ stabilizes in the sense of the definition of stabilization given in [7], then it is also stabilizing in the sense of the definition given here, by means of an estimate of the type (9). (The converse implication is obvious.) Suppose $M(R) > 0$, $\delta_0(r, R) > 0$, and $T(r, R) \geq 0$, are such that $\lim_{R \searrow 0} M(R) = 0$, and whenever $\mathbf{d}(\pi) < \delta_0(r, R)$ and $|x^0| \leq R$, necessarily $|x_\pi(t, x^0)| \leq M(R)$ for all $t \geq 0$ and $|x_\pi(t, x^0)| \leq r$ for all $t \geq T(r, R)$. We first define $\delta = \delta(\varepsilon, R) > 0$, for each $0 < \varepsilon < R$, as follows: pick the smallest positive integer $k$ (necessarily $\geq 2$) such that $\frac{R}{k} \leq \varepsilon$, and let

$$\delta(\varepsilon, R) := \min\left\{\delta_0\left(\frac{R}{2}, R\right), \ldots, \delta_0\left(\frac{R}{k}, R\right)\right\}.$$

Next, we define a function

$$\varphi : \mathbb{R}_{\geq 0} \times \mathbb{R}_{\geq 0} \to \mathbb{R}_{\geq 0}$$

as follows. Pick any $R > 0$ and any $t \geq 0$. Let $0 < t_1 < t_2 < \ldots$ be a sequence of real numbers (depending on $R$) so that $t_i \to \infty$ as $i \to \infty$ and such that

$$t_i \geq T\left(\frac{R}{i+1}, R\right), \quad i = 1, 2, \ldots.$$

Now define $\varphi(R, t) := M(R)$ for $t \in [0, t_1)$ and $\varphi(R, t) := \frac{R}{i+1}$ if $t \in [t_i, t_{i+1})$ for some $i \geq 1$.

We claim now: for each $0 < \varepsilon < R$, each $|x^0| \leq R$, each $t \geq 0$, and each sampling schedule such that $\mathbf{d}(\pi) \leq \delta(\varepsilon, R)$,

$$|x_\pi(t, x^0)| \leq \max\{\varphi(R, t), \varepsilon\}.$$

Pick any such $\varepsilon, R, x^0, \pi$. Define the sequence $\{t_i\}$ as above, for this $R$. Let $k$ be the smallest positive integer such that $\frac{R}{k} \leq \varepsilon$. By definition of $\delta$, $\delta(\varepsilon, R) \leq \delta_0(R/j, R)$ for all $j = 2, \ldots, k$. We consider three cases: (i) $t < t_1$, (ii) $t \in [t_1, t_{k-1})$, and (iii) $t \geq t_{k-1}$. In the first case, we know that $|x_\pi(t, x^0)| \leq M(R) = \varphi(R, t)$. In the last case, we know that $|x_\pi(t, x^0)| \leq \varepsilon$, because $\mathbf{d}(\pi) \leq \delta_0(R/k, R)$ implies that $|x_\pi(t, x^0)| \leq R/k \leq \varepsilon$ for all $t \geq T(R/k, R)$, and $t_{k-1} \geq T(R/k, R)$ by definition of the $t_i$'s. In case (ii), we have that there is some $j \in \{2, \ldots, k-1\}$ so that $t \in [t_{j-1}, t_j)$. Since $t_{j-1} \geq T(R/j, R)$ and $\mathbf{d}(\pi) \leq \delta_0(R/j, R)$, $|x_\pi(t, x^0)| \leq R/j = \varphi(R, t)$. The claim is then established.

It only remains to show that there is some function $\beta \in \mathcal{KL}$ so that $\varphi(s, t) \leq \beta(s, t)$ for all $s, t$. The constructions given in the first section of [28] show that there exists such a $\beta$, provided that the following properties hold for $\varphi$:

1. For some $\gamma \in \mathcal{K}_\infty$, $\sup_{t>0} \varphi(R, t) \leq \gamma(R)$ for all $R > 0$.

2. For each $\varepsilon > 0$ and each $R > 0$ there is some $T(\varepsilon, R)$ such that

$$t \geq T(\varepsilon, R) \;\Rightarrow\; \varphi(R, t) \leq \varepsilon.$$

The first property is satisfied because $M$ could, without loss of generality, be taken to be of class $\mathcal{K}_\infty$ (if necessary, first use a step function $M$; then majorize it by a $\mathcal{K}_\infty$ function). The second property holds as well: take without loss of generality $\varepsilon < R$, then pick a positive integer $k$ minimal with $R/k \leq \varepsilon$, let the sequence $\{t_i\}$ be as in the construction of $\varphi$, and define $T(\varepsilon, R) := t_{k-1}$. Observe that $t \geq t_{k-1}$ implies, by definition of $\varphi$, that $\varphi(R,t) \leq R/k \leq \varepsilon$, as wanted.

## 5.2 Proximal Form of clf Theorem

We fill-in the details here to show how the proximal subgradient form of the continuous clf existence Theorem D follows from the result in [34]. Before stating the original form of the result, we recall the notion of relaxed control. For each real $s > 0$, let us denote by $\mathcal{U}_s$ the radius-$s$ ball in $\mathcal{U} = \mathbb{R}^m$. A *relaxed $\mathcal{U}_s$-valued control* is a measurable map $\omega : I \to \mathbb{P}(\mathcal{U}_s)$, where $I$ is an interval containing zero and $\mathbb{P}(\mathcal{U}_s)$ denotes the set of all Borel probability measures on $\mathcal{U}_s$. Note that ordinary controls can be seen also as relaxed controls, via the natural embedding of $\mathcal{U}_s$ into $\mathbb{P}(\mathcal{U}_s)$ (map any point $u \in \mathcal{U}_s$ into the Dirac delta measure supported at $u$). Given any $\mu \in \mathbb{P}(\mathcal{U}_s)$, we write $\int_{\mathcal{U}_s} f(x,u)\,d\mu(u)$ simply as $f(x,\mu)$. As with ordinary controls, we also denote by $x(t, x^0, \omega)$ the solution of the initial value problem that obtains from initial state $x^0$ and relaxed control $\omega$, and we consider the supremum norm $\|\omega\|$, defined as the infimum of the set of $s$ such that $\omega(t) \in \mathbb{P}(\mathcal{U}_s)$ for almost all $t \in I$.

The main result in [34] says that if (and only if) a system is AC, there exist two continuous, positive definite functions $V, W : \mathbb{R}^n \to \mathbb{R}$, with $V$ proper, and a nondecreasing function $\sigma : \mathbb{R}_{\geq 0} \to \mathbb{R}_{\geq 0}$, so that the following property holds: for each $x^0 \in \mathbb{R}^n$ there are a $T > 0$ and a relaxed control $\omega : [0, T) \to \mathbb{P}(\mathcal{U}_{\sigma(|x^0|)})$, so that $x(t) := x(t, x^0, \omega)$ is defined for all $0 \leq t < T$ and

$$V(x(t)) - V(x^0) \leq -\int_0^t W(x(\tau))\,d\tau \quad \forall t \in [0, T). \tag{16}$$

To simplify notations, let us write $s := \sigma(|x^0|)$. In order to obtain the proximal version of the result, we show that

$$\min_{|u| \leq s} \zeta \cdot f(x^0, u) \leq -W(x^0) \tag{17}$$

for all $\zeta \in \partial_P V(x^0)$.

As in [40], we first make an intermediate reduction, showing that

$$\min_{v \in F(x^0, s)} DV(x^0; v) \leq -W(x^0), \tag{18}$$

where $F(x^0, s)$ denotes the (closed) convex hull of $\{f(x^0, u), u \in \mathcal{U}_s\}$, and $DV(x^0; v)$ is the directional subderivate, or contingent epiderivative, of $V$ in the direction of $v$ at $x^0$, defined as

$$DV(x^0; v) := \liminf_{\substack{t \searrow 0 \\ v' \to v}} \frac{1}{t}[V(x^0 + tv') - V(x^0)].$$

(The minimum in (18) is achieved, because the map $v \mapsto DV(x^0; v)$ is lower semicontinuous, see e.g. [6], ex.3.4.1e.) So, let $x(t) = x(t, x^0, \omega)$ be as above, and consider for each $t \in [0, T)$ the vectors

$$r_t := \frac{1}{t}(x(t) - x^0) = \frac{1}{t}\int_0^t f(x(\tau), \omega(\tau))\,d\tau = q_t + p_t,$$

with
$$p_t := \frac{1}{t} \int_0^t f(x^0, \omega(\tau)) \, d\tau \, ,$$

where $q_t \to 0$ as $t \searrow 0$ (the existence of such a $q_t$, for each $t$, is an easy consequence of the fact that $f$ is locally Lipschitz on $x$, uniformly on $u \in \mathcal{U}_s$). Moreover, $p_t \in F(x^0, s)$ for all $t$ (because $F(x^0, s)$ is convex, so $f(x^0, \omega(\tau)) \in F(x^0, s)$ for each $\tau$, and then using convexity once more). By compactness, we have that there is some $v \in F(x^0, s)$ and some subsequence $p_{t_j} \to v$ with $t_j \searrow 0$; as $q_t \to 0$, also $v_j := r_{t_j} \to v$. For this $v$,

$$\begin{aligned} DV(x^0; v) &\leq \liminf_{j \to \infty} \frac{1}{t_j}[V(x^0 + t_j v_j) - V(x^0)] \\ &= \liminf_{j \to \infty} \frac{1}{t_j}[V(x(t_j)) - V(x^0)] \leq -W(x^0) \end{aligned}$$

(using (16)), so (18) indeed holds.

Finally, we show that (18) implies (17). Let $v \in F(x^0, s)$ achieve the minimum, and pick any $\zeta \in \partial_\mathrm{P} V(x^0)$. By definition of proximal subgradient, there is some $\mu > 0$ so that, for each $v' \in F(x^0, s)$ and each small $t \geq 0$,

$$\zeta \cdot v' \leq \frac{1}{t}[V(x^0 + tv') - V(x^0)] + \mu t \, |v'|^2$$

so taking limits along any sequence $t \searrow 0$ and $v' \to v$ shows that $\zeta \cdot v \leq DV(x^0; v) \leq -W(x^0)$. By definition of $F(x^0, s)$, this means that there must exist a $u \in \mathcal{U}_0$ so that also $\zeta \cdot f(x^0, u) \leq -W(x^0)$, as desired.

Observe that the proximal condition (17), being linear on the velocities $f(x, u)$, has the great advantage of not requiring convex hulls in order to state, and in that sense is far more elegant than (18).

### 5.3 $V_\alpha$ is a (Local) clf

We prove here that, on each set $C_{r,R} = \{x \in \mathbb{R}^n \mid r \leq |x| \leq R\}$, there is an $\alpha_0(r, R) > 0$ such that, for each positive $\alpha \leq \alpha_0(r, R)$, the (locally Lipschitz) function $V_\alpha$ behaves like a clf on the set $C_{r,R}$. In order to simplify referencing, we write "([7]n)" to refer to Equation (n) in [7] and do not redefine notations given there. Let $\mathcal{U}_0$ be a compact subset so that $\min_{u \in \mathcal{U}_0} \zeta \cdot f(x, u) \leq -W(x)$ for every $\zeta \in \partial_\mathrm{P} V(x)$ and every $x$ in the ball of radius $R + \sqrt{2\beta(R)}$. Let $m_{r,R} := \frac{1}{2} \min\{W(x) \mid x \in C_{r,R}\}$, and let $\ell$ be so that $|f(x, u) - f(x', u)| \leq \ell |x - x'|$ for all $u \in \mathcal{U}_0$ and all $x$ in the ball of radius $R + \sqrt{2\beta(R)}$. (This is almost as in ([7]29), except that, there, $\ell$ was a Lipschitz constant only with respect to $|x| \leq R$.) Finally, as in [7], $\omega_R(\cdot)$ denotes the modulus of continuity of $V$ on the ball of radius $R + \sqrt{2\beta(R)}$.

**Proposition.** Let $\alpha \in (0, 1]$ satisfy $2\ell \omega_R\left(\alpha\sqrt{2\beta(R)}\right) \leq m_{r,R}$. Then, for all $x \in C_{r,R}$ and all $\zeta \in \partial_\mathrm{P} V_\alpha(x)$,
$$\min_{u \in \mathcal{U}_0} \zeta \cdot f(x, u) \leq -m_{r,R} \, .$$

*Proof.* We start by remarking that, for all $\alpha \in (0, 1]$ and all $x \in B_R$,

$$|y_\alpha(x) - x|^2 \leq 2\alpha^2 \, \omega_R\left(\alpha\sqrt{2\beta(R)}\right) \, . \tag{19}$$

This follows from (cf. ([7]19)):

$$\frac{1}{2\alpha^2} |y_\alpha(x) - x|^2 \leq V(x) - V(y_\alpha(x)) \leq \omega_R(|y_\alpha(x) - x|)$$

and using that $|y_\alpha(x) - x| \leq \alpha \sqrt{2\beta(R)}$ (cf. Lemma III.3 in [7]) plus the fact that $\omega_R(\cdot)$ is nondecreasing. So

$$|\zeta_\alpha(x)| \cdot |y_\alpha(x) - x| = \frac{|y_\alpha(x) - x|^2}{\alpha^2} \leq 2\omega_R\left(\alpha\sqrt{2\beta(R)}\right)$$

and hence

$$\zeta_\alpha(x) \cdot [f(y_\alpha(x), u) - f(x, u)] \leq \ell |\zeta_\alpha(x)| \cdot |y_\alpha(x) - x| \leq 2\ell \omega_R\left(\alpha\sqrt{2\beta(R)}\right)$$

for all $u \in \mathcal{U}_0$ and all $x \in B_R$, because $y_\alpha(x) \in B_{R+\sqrt{2\beta(R)}}$ (cf. ([7]20)). Thus,

$$\zeta_\alpha(x) \cdot f(x, u) \leq \zeta_\alpha(x) \cdot f(y_\alpha(x), u) + 2\ell \omega_R\left(\alpha\sqrt{2\beta(R)}\right).$$

For each $x \in C_{r,R}$, we have that $y_\alpha(x) \in B_{R+\sqrt{2\beta(R)}}$ and (cf. ([7]12)) $\zeta_\alpha(x) \in \partial_P V(y_\alpha(x))$, so we can pick a $u_x \in \mathcal{U}_0$ so that $\zeta_\alpha(x) \cdot f(y_\alpha(x), u_x) \leq -W(x) \leq -2m_{r,R}$. Using now the fact that $2\ell \omega_R\left(\alpha\sqrt{2\beta(R)}\right) \leq m_{r,R}$, we conclude that

$$\zeta_\alpha(x) \cdot f(x, u) \leq -m_{r,R}.$$

So we only need to prove that $\partial_P V_\alpha(x) \subseteq \{\zeta_\alpha(x)\}$ for all $x$. This is given in [6], Theorem 1.'5.1, but is easy to prove directly: pick any $\zeta \in \partial_P V_\alpha(x)$; by definition, this means that exists some $\sigma > 0$ such that, for all $y$ near $x$, $\zeta \cdot (y - x) \leq V_\alpha(y) - V_\alpha(x) + \sigma |y - x|^2$, so

$$\zeta \cdot (y - x) \leq V_\alpha(y) - V_\alpha(x) + \gamma(|y - x|) \qquad (20)$$

for some $\gamma(r) = o(r)$. Adding

$$-\zeta_\alpha(x) \cdot (y - x) \leq -V_\alpha(y) + V_\alpha(x) + \frac{1}{2\alpha^2} |y - x|^2$$

(cf. ([7]13)), we conclude

$$(\zeta - \zeta_\alpha(x)) \cdot (y - x) \leq o(|y - x|).$$

Substituting $y = x + h(\zeta - \zeta_\alpha(x))$ and letting $h \searrow 0$ shows that $\zeta = \zeta_\alpha(x)$. (Observe that we have proved more that claimed: Equation (20) is satisfied by any viscosity subgradient $\zeta$; in particular, the gradient of $V$, if $V$ happens to be differentiable at the point $x$, must coincide with $\zeta_\alpha(x)$.) ∎

Observe that the feedback constructed in [7] was $k(x) =$ any $u \in \mathcal{U}_0$ minimizing $\zeta_\alpha(x) \cdot f(x, u)$ (where $\alpha$ and $\mathcal{U}_0$ are chosen constant on certain compacts). As $V_\alpha$ is locally Lipschitz, it is differentiable almost everywhere. Thus, the Proposition (see the end of the proof) insures that $\zeta_\alpha(x) = \nabla V_\alpha(x)$ for almost all $x$. So $k(x) = u$ is, at those points, the pointwise minimizer of the Hamiltonian $\nabla V_\alpha(x) \cdot f(x, u)$ associated to the regularized clf $V_\alpha$.

## 5.4 Proof of Theorem E

We will prove the following more precise result. All undefined notations, including the definitions of the functions $\gamma$ and $\rho$, can be found in [7].

**Theorem G** *Pick any $0 < r < R$ so that $2\gamma(r) < \gamma(R)$. Then there exist positive numbers $\alpha, \delta, \kappa, T$ such that, for each partition $\pi$ with $\mathbf{d}(\pi) \leq \delta$, and each $e : [0, \infty) \to \mathbb{R}^n$ which satisfies $|e(t)| \leq \kappa \underline{\mathbf{d}}(\pi)$ for all $t$, the following property holds: if $x(\cdot)$ satisfies*

$$\dot{x} = f(x, k(x+e)), \quad |x(0)| \leq \tfrac{1}{2}\rho(R), \tag{21}$$

*where $k$ is the feedback $k_{\alpha,r,R}$, then*

$$x(t) \in B_R \quad \forall t \geq 0 \tag{22}$$

*and*

$$x(t) \in B_r \quad \forall t \geq T. \tag{23}$$

Theorem E is a corollary: we first pick any $\Gamma$ such that $s \leq (1/2)\rho(\Gamma(s))$ for all $s \geq 0$. Then, given $\varepsilon$ and $K$, we can let $R := \Gamma(K)$ (so, $K \leq (1/2)\rho(R)$), we then take any $0 < r < \varepsilon$ such that $2\gamma(r) \leq \gamma(R)$, and apply the above result to $0 < r < R$.

We prove Theorem G through a series of technical steps. Let $0 < r < R$ be given, with $2\gamma(r) < \gamma(R)$. In order to simplify referencing, we write "([7]n)" to refer to Equation (n) in [7].

We start by picking $\alpha$ as any positive number which satisfies ([7]23), ([7]30), and, also, instead of ([7]40), the slightly stronger condition

$$\omega_R\left(\sqrt{2\beta(R)}\alpha\right) < \frac{1}{16}\gamma(r). \tag{24}$$

The function $V_\alpha$ is defined as in ([7]9); because of ([7]13), ([7]11), and ([7]19), the number

$$c := \frac{\sqrt{2\beta(R)}}{\alpha} + \frac{R}{\alpha^2}$$

is a Lipschitz constant for $V_\alpha$ on the set $B_R$. Without loss of generality, we assume $c \geq 2$. We let $\mathcal{U}_0$ be as in [7], and take the feedback $k = k_{\alpha,r,R}$. The numbers $\ell$, $m$, and $\Delta$ are as in ([7]29) and the equation that follows it.

Next, we pick any $\varepsilon_0 > 0$ so that all the following – somewhat redundant – properties hold:

$$B_{\frac{1}{2}\rho(R)} + \varepsilon_0 B \subseteq G_R^\alpha \tag{25}$$

(this is possible because $B_{\rho(R)} \subseteq G_R^\alpha$ by ([7]22)),

$$G_R^\alpha + 2\varepsilon_0 B \subseteq B_R \tag{26}$$

(possible because $G_R^\alpha \subseteq \text{int } B_R$ by ([7]24)),

$$\varepsilon_0 \leq \frac{\gamma(r)}{8c} \leq \frac{\gamma(r)}{16}, \tag{27}$$

$$\gamma(r) + c\,\varepsilon_0 \;<\; \tfrac{1}{2}\gamma(R) \tag{28}$$

(recall that $2\gamma(r) < \gamma(R)$), and

$$G_r^\alpha + 2\varepsilon_0 B \;\subseteq\; G_R^\alpha. \tag{29}$$

We let $\delta_0$ be any positive number so that, for every initial state $x^o$ in the compact set $G_R^\alpha$, and for each control $u : [0,\delta_0] \to \mathcal{U}_0$, the solution of $\dot x = f(x,u)$ with $x(0) = x^o$ is defined on the entire interval $[0,\delta_0]$ and satisfies $x(t) \in G_R^\alpha + \varepsilon_0 B$ for all $t$.

Finally, we pick any $\delta > 0$ which satisfies ([7]33), ([7]41), as well as

$$\delta \;\leq\; \min\left\{1,\, \delta_0,\, \frac{2}{\Delta}\varepsilon_0,\, \frac{\gamma(r)}{8cm}\right\}, \tag{30}$$

and we let

$$\kappa \;:=\; \frac{\Delta}{4c} e^{-\ell} \tag{31}$$

and

$$T \;:=\; \frac{\gamma(R)}{\Delta} \tag{32}$$

(this is twice the value used in ([7]39)).

The main technical fact needed is as follows.

**Proposition 5.4.1** *Let $0 < r < R$ satisfy $2\gamma(r) < \gamma(R)$. Let $\alpha, \delta, \kappa, T$ be defined as above, and $k$ be the feedback $k_{\alpha,r,R}$. Pick any $\varepsilon > 0$, and consider the following set:*

$$P = P_{r,R,\varepsilon} \;:=\; \{x \mid x + \varepsilon B \subseteq G_R^\alpha\}\,.$$

*Let $\pi$ be a partition which satisfies*

$$\frac{\varepsilon}{\kappa} \;\leq\; t_{i+1} - t_i \;\leq\; \delta \quad \forall\, i = 0, 1, \ldots \tag{33}$$

*(that is, $\mathbf{d}(\pi) \leq \delta$ and $\varepsilon \leq \kappa\,\underline{\mathbf{d}}(\pi)$). Then, for any $e : [0,\infty) \to \mathbb{R}^n$ such that $|e(t)| \leq \varepsilon$ for all $t$, and any $x^o \in P$, the solution of $\dot x = f(x, k(x+e))$ is defined for all $t \geq 0$, it satisfies (22) and (23), and $x(t_i) \in P$ for all $i$.*

We will prove this via a couple of lemmas, but let us first point out how Theorem G follows from Proposition 5.4.1. Suppose given a partition $\pi$ with $\mathbf{d}(\pi) \leq \delta$, an error function $e$ so that $|e(t)| \leq \kappa\,\underline{\mathbf{d}}(\pi)$ for all $t$, and an initial state $x^o$ with $|x^o| \leq \tfrac{1}{2}\rho(R)$. Let $\varepsilon := \sup_{t \geq 0} |e(t)| \leq \kappa\,\underline{\mathbf{d}}(\pi)$. Note that

$$\varepsilon \;\leq\; \kappa\,\underline{\mathbf{d}}(\pi) \;\leq\; \kappa\,\mathbf{d}(\pi) \;\leq\; \kappa\,\delta \;\leq\; \kappa\frac{2}{\Delta}\varepsilon_0 \;\leq\; \kappa\frac{4c}{\Delta}e^{\ell}\varepsilon_0 \;=\; \varepsilon_0\,. \tag{34}$$

Thus, by (25),

$$B_{\frac{1}{2}\rho(R)} \;\subseteq\; P\,, \tag{35}$$

so $x^o \in P$. Therefore, (22) and (23) hold for the solution of $\dot x = f(x, k(x+e))$, as wanted for Theorem G.

We now prove Proposition 5.4.1. Observe that (33) implies $\varepsilon \leq \kappa\,\underline{\mathbf{d}}(\pi)$, so, arguing as in (34), $\varepsilon \leq \varepsilon_0$. It is useful to introduce the following set as well:

$$Q = Q_{r,R,\varepsilon} \;:=\; G_r^\alpha + \varepsilon B\,.$$

Observe that $Q \subseteq P$, by (29).

We start the proof by establishing an analogue of Lemma IV.2 in [7]:

**Lemma 5.4.1** *If, for some index $i$, $x_i := x(t_i) \in P \setminus Q$, then $x(t)$ is defined for all $t \in [t_i, t_{i+1}]$,*

$$x(t) \in B_R \quad \forall t \in [t_i, t_{i+1}], \tag{36}$$

$$V_\alpha(x(t)) \leq V_\alpha(x_i) + \varepsilon_0 \quad \forall t \in [t_i, t_{i+1}], \tag{37}$$

*and, letting $x_{i+1} := x(t_{i+1})$:*

$$x_{i+1} \in P, \tag{38}$$

$$V_\alpha(x_{i+1}) - V_\alpha(x_i) \leq -\frac{\Delta}{2}(t_{i+1} - t_i). \tag{39}$$

*Proof.* By definition of $P$, we have that $\tilde{x}_i := x_i + e(t_i) \in G_R^\alpha$. Also, $\tilde{x}_i \notin G_r^\alpha$, since otherwise $x_i$ would belong to $Q$. In particular, $\tilde{x}_i \notin B_{\rho(r)}$. Let $\tilde{x}(\cdot)$ be the solution of $\dot{x} = f(x, k(\tilde{x}_i))$ with $\tilde{x}(t_i) = \tilde{x}_i$ on $[t_i, t_{i+1}]$. By Lemma IV.2 in [7], this solution is well-defined and it holds that $\tilde{x}(t_{i+1}) \in G_R^\alpha$ and $V_\alpha(\tilde{x}(t)) - V_\alpha(\tilde{x}_i) \leq -\Delta(t - t_i)$ for all $t$. As $x_i \in P \subseteq G_R^\alpha$, and $t_{i+1} - t_i \leq \mathbf{d}(\pi) \leq \delta \leq \delta_0$, the definition of $\delta_0$ insures that the solution $x(\cdot)$ of $\dot{x} = f(x, k(\tilde{x}_i))$ with $x(t_i) = x_i$ is indeed well-defined, and it stays in $G_R^\alpha + \varepsilon_0 B \subseteq B_R$ for all $t$. So, by Gronwall's inequality, we know that

$$|x(t) - \tilde{x}(t)| \leq e^{(t-t_i)\ell}|x_i - \tilde{x}_i| \leq e^{\delta\ell}\varepsilon$$

for all $t \in [t_i, t_{i+1}]$. Since $V_\alpha$ has Lipschitz constant $c$ on $B_R$, we have

$$\begin{aligned}
V_\alpha(x(t)) - V_\alpha(x_i) &\\
&= V_\alpha(x(t)) - V_\alpha(\tilde{x}(t)) + V_\alpha(\tilde{x}(t)) - V_\alpha(\tilde{x}_i) + V_\alpha(\tilde{x}_i) - V_\alpha(x_i)\\
&\leq ce^{\delta\ell}\varepsilon - \Delta(t - t_i) + c\varepsilon\\
&\leq \frac{\Delta}{2}(t_{i+1} - t_i) - \Delta(t - t_i)
\end{aligned}$$

for all $t \in [t_i, t_{i+1}]$, where we have used that $\delta \leq 1$, $\varepsilon \leq \kappa\underline{\mathbf{d}}(\pi)$, the definition $\kappa = (\Delta/4c)e^{-\ell}$, and the fact that $\underline{\mathbf{d}}(\pi) \leq t_{i+1} - t_i$. In particular, the estimate (39) results at $t = t_{i+1}$, and (37) holds because $(\Delta/2)\delta \leq \varepsilon_0$ by (30).

We are only left to prove that $x_{i+1} \in P$. By definition of $P$, this means that for any given $\eta \in \mathbb{R}^n$ with $|\eta| \leq \varepsilon$ it must hold that

$$\eta + x_{i+1} \in G_R^\alpha = \{x \mid V_\alpha(x) \leq (1/2)\gamma(R)\}.$$

Pick such an $\eta$; then $|\eta + x_{i+1} - \tilde{x}(t_{i+1})| \leq \varepsilon + e^{\delta\ell}\varepsilon \leq 2e^\ell\varepsilon$, and so, since $\eta + x_{i+1} \in G_R^\alpha + 2\varepsilon_0 B \subseteq B_R$,

$$\begin{aligned}
V_\alpha(\eta + x_{i+1}) &\leq 2ce^\ell\varepsilon + V_\alpha(\tilde{x}(t_{i+1}))\\
&\leq 2ce^\ell\varepsilon + V_\alpha(\tilde{x}_i) - \Delta(t_{i+1} - t_i)\\
&\leq V_\alpha(\tilde{x}_i) - \frac{\Delta}{2}(t_{i+1} - t_i) \leq \frac{1}{2}\gamma(R)
\end{aligned}$$

where we again used $\varepsilon \leq \kappa\underline{\mathbf{d}}(\pi)$ as well as the fact that $V_\alpha(\tilde{x}_i) \leq (1/2)\gamma(R)$ (since $\tilde{x}_i \in G_R^\alpha$). ∎

We also need another observation, this one paralleling the proof of Lemma IV.4 in [7].

**Lemma 5.4.2** *If, for some index $i$, $x_i := x(t_i) \in Q$, then $x(t)$ is defined for all $t \in [t_i, t_{i+1}]$,*

$$V_\alpha(x(t)) \leq \frac{3}{4}\gamma(r) \quad \forall t \in [t_i, t_{i+1}], \tag{40}$$

*and*

$$V(x(t)) \leq \frac{7}{8}\gamma(r) \quad \forall t \in [t_i, t_{i+1}]. \tag{41}$$

*In particular, $x(t) \in B_r$ for all $t \in [t_i, t_{i+1}]$ and $x(t_{i+1}) \in P$.*

*Proof.* The fact that $x$ is defined follows from the choice of $\delta_0$, and we know that $x(t) \in B_R$ for all $t \in [t_i, t_{i+1}]$. So

$$|x(t) - x_i| \leq m\delta \quad \forall t \in [t_i, t_{i+1}]. \tag{42}$$

By definition of $Q$, we may write $x_i = x' + \eta$, for some $x' \in G_r^\alpha$ and some $|\eta| \leq \varepsilon$. Thus, $V_\alpha(x_i) \leq V_\alpha(x') + c\varepsilon_0 \leq \frac{1}{2}\gamma(r) + c\varepsilon_0$ (second inequality by definition of $G_r^\alpha$). Together with (42), this gives

$$V_\alpha(x(t)) \leq \frac{1}{2}\gamma(r) + c\varepsilon_0 + cm\delta \leq \frac{3}{4}\gamma(r) \quad \forall t \in [t_i, t_{i+1}]$$

(using (27) and (30)). So, using (24) and ([7]21),

$$V(x(t)) \leq V_\alpha(x(t)) + \omega_R\left(\sqrt{2\beta(R)}\alpha\right) \leq \frac{3}{4}\gamma(r) + \frac{1}{16}\gamma(r) < \frac{7}{8}\gamma(r)$$

for all $t \in [t_i, t_{i+1}]$, as wanted. By the definition of $\gamma$ (see [7]), this means that $x(t) \in B_r$ for all $t \in [t_i, t_{i+1}]$. Finally, if $|\eta| \leq \varepsilon$ then

$$V_\alpha(x(t_{i+1}) + \eta) \leq c\varepsilon + \frac{3}{4}\gamma(r) \leq \frac{1}{2}\gamma(R)$$

(the last inequality by (28)), which means that $x(t_{i+1}) + \eta \in G_R^\alpha$; this implies that $x(t_{i+1}) \in P$. ∎

Back to the proof of Proposition 5.4.1, since $x^0 \in P$, Lemmas 5.4.1 and 5.4.2 guarantee that the solution exists for all $t$ and remains in $B_R$, and that $x_i := x(t_i) \in P$ for all $i$.

Moreover, if there is some $j$ so that $x(t_j) \in Q$, then it holds that $V_\alpha(x_i) \leq \frac{3}{4}\gamma(r)$ for all $i > j$. This is because on intervals in which $x_{i-1} \in Q$, we already know that $V_\alpha(x(t_i)) \leq \frac{3}{4}\gamma(r)$, and if instead $x_{i-1} \in P \setminus Q$, then we have $V_\alpha(x_i) < V_\alpha(x_{i-1})$. So, for any such $i > j$, either $x(t) \in B_r$ for all $t \in [t_i, t_{i+1}]$ (first case) or $V_\alpha(x(t)) \leq \frac{3}{4}\gamma(r) + \varepsilon_0$ for all $t \in [t_i, t_{i+1}]$ (second case). Actually, in this last case we also have

$$\begin{aligned} V(x(t)) &\leq V_\alpha(x(t)) + \omega_R\left(\sqrt{2\beta(R)}\alpha\right) \\ &\leq \frac{3}{4}\gamma(r) + \varepsilon_0 + \omega_R\left(\sqrt{2\beta(R)}\alpha\right) \\ &\leq \frac{3}{4}\gamma(r) + \frac{1}{16}\gamma(r) + \frac{1}{16}\gamma(r) < \gamma(r), \end{aligned}$$

(using (27) and (24)), so, again by definition of $\gamma$, also $x(t) \in B_r$ for all $t \in [t_i, t_{i+1}]$. In conclusion, trajectories stay in $B_r$ after the first time that $x(t_j) \in Q$. So we only need to show that there is such a $j$, with $t_j \leq T$.

Suppose instead that for $i = 0, \ldots, k$ it holds that $x(t_i) \notin Q$, and $t_k > T$. Applying (39) repeatedly,

$$0 \leq V_\alpha(x(t_k)) \leq V_\alpha(x^0) - \frac{\Delta}{2} t_k < \frac{\gamma(R)}{2} - \frac{\Delta}{2} T,$$

(recall that $x^0 \in P$ implies that $x \in G_R^\alpha$), and this contradicts (32). The proof of Proposition 5.4.1 is then complete. ∎

This completes the proof. We remark that a more global result is also possible, as follows. We start by picking two sequences $\{r_j, j \in \mathbb{Z}\}$ and $\{R_j, j \in \mathbb{Z}\}$ such that $r_j, R_j \to 0$ as $j \to -\infty$, $r_j, R_j \to \infty$ as $j \to \infty$, and $2R_j \leq \rho(R_{j+1})$, $2\gamma(r_j) < \gamma(R_j)$, and $2r_j \leq \rho(R_{j-1})$ for all $j$. Next we pick, for each $j$, positive numbers $\alpha_j, \delta_j, \kappa_j, T_j$ associated as per Proposition 5.4.1 to $r_j$ and $R_j$, and let $k_j := k_{\alpha_j, r_j, R_j}$. (We may assume that values of $k_j$ belong to some fixed $\mathcal{U}_0$ for all $j \leq 0$, and to $\mathcal{U}_j$ for $j > 0$, with all the $\mathcal{U}_j$ compact and forming an increasing sequence.) Since $G_{R_j}^{\alpha_j} \subseteq \mathrm{int}\, G_{R_{j+1}}^{\alpha_{j+1}}$ for all $j$ (this is proved just as in ([7]48)), there is some sequence of positive numbers $\{\varepsilon_j, j \in \mathbb{Z}\}$ so that $\varepsilon_j < \kappa_j \delta_j$ for all $j$ and also, denoting $P_j := P_{r_j, R_j, \varepsilon_j}$, so that $P_j \subseteq P_{j+1}$ for all $j$. Note that the choice of the $r_j$'s and $R_j$'s assures that $\bigcup P_j = \mathbb{R}^n \setminus \{0\}$. Since $2r_j \leq \rho(R_{j-1})$ for all $j$, and using (35), we know that $B_{r_j} \subseteq P_{j-1}$ for all $j$.

Finally, we define the feedback $k: \mathbb{R}^n \to \mathbb{R}^m$ via:

$$k(x) := k_j(x) \text{ if } x \in P_j \setminus P_{j-1}$$

(and $k(0) = 0$), and let $\varepsilon(x) := \varepsilon_j$, $\mu(x) := \frac{\varepsilon_j}{\kappa_j}$, and $\delta(x) := \delta_j$ for $x \in P_j \setminus P_{j-1}$, for each $j$ (let $\varepsilon_j(0) = 0$ and $0 < m(0) < \delta(0)$ be arbitrary). Observe that, since $\varepsilon_j < \kappa_j \delta_j$ for all $j$, it holds that $\mu(x) < \delta(x)$ for all $x$.

Now suppose that $x^0$ is given, and a partition $\pi$ and a function $e(\cdot)$ are given so that, recursively along the solution of $\dot{x} = f(x, k(x+e))$,

$$\mu(x(t_i)) \leq t_{i+1} - t_i \leq \delta(x(t_i))$$

and

$$e(t_i) \leq \varepsilon(x(t_i))$$

for all $i$. Note that, if for some $i$ and $j$, $x(t_i) \in P_j \setminus P_{j-1}$, then the inequality in (33) holds with $\varepsilon = \varepsilon_j$, $\kappa = \kappa_j$, $\delta = \delta_j$, and this $i$.

Let $j \in \mathbb{Z}$ be maximal so that $x^0 \in P_j$. Then, it will hold that the solution stays in the bounded set $B_{R_j}$, and every $x(t_i)$ is again in $P_j$, until the first sampling time $t_i$ that $x(t_i) \in P_\ell$ for some $\ell < j$. This first time, say $t_q$, is at most $T_j$, because $B_{r_j} \subseteq P_{j-1}$. If $x(t_q) = 0$, then the solution stays there forever (since $f(0,0) = 0$). Otherwise, $x(t_q) \in P_\ell \setminus P_{\ell-1}$ for some $\ell < j$, and we may repeat the argument. The conclusion is that the trajectory keeps visiting smaller sets $P_\ell$ (or it becomes 0 in finite time), with an upper bound $T(i,j)$ (the sum of the corresponding $T_\ell$'s) on the time required for entering a given $P_i$, if the initial state $x^0$ was in a given $P_j$. Furthermore, given any $0 < r < R$, there are $i < j$ so that $P_i \subseteq B_r \subseteq B_R \subseteq P_j$. Thus, all trajectories in $B_R$ are taken into $B_r$ in a uniform time $T(r,R)$, with bounded overshoot (since trajectories stay in $B_{R_j}$).

## 6  Stability with Respect to Disturbances

One of the conclusions from the preceding discussion is that, provided that there is a $\mathcal{C}^1$ clf, the closed-loop system $\dot{x} = f(x, k(x+e)+d)$ in Equation (13) is "stable" in a certain precisely defined sense. If we view the measurement and actuator error functions $e$ and $d$ as external inputs to the system, we may write $\dot{x} = g(x, v)$ to represent this closed-loop system (13), where $v = (e, d)$. Interpreted in these terms, the previous discussion dealt with the effect of *small* "inputs" $v$ on the system $\dot{x} = g(x, v)$.

In many control problems, on the other hand, one wishes to understand the effect of *large* external inputs. This interest might be due to the need to quantity the effect of large errors in measurements or actuators, but there are many other motivations as well. For example, in disturbance rejection problems, one may wish to cancel the effect of an input to the system (such as wind acting on an aircraft). In tracking problems, a desired path $r(t)$ may be given and we might want to make the tracking error $x(t) - r(t) \to 0$ as $t \to \infty$, which can be interpreted loosely as the stabilization of the signal $x - r$. (More generally, we may want a subset of the coordinates, or even an arbitrary function of the state, to follow a desired behavior; see [36] for a discussion of some tracking problems.) We are thus led to the study of the general question: how does the state trajectory $x(t) = x(t, x^o, u)$ depend on $u$ and $x^o$? More precisely, we would like to compare the global behavior of $x$ and $u$, asking questions such as: does boundedness of $u$ imply boundedness of $x$?, or does $u(t) \to 0$ as $t \to \infty$ imply $x(t) \to 0$ as $t \to \infty$? And if yes, how does the initial state affect the bounds or convergence rates, respectively?

There are two very conceptually different ways of formulating the notion of stability of control systems. One of them, which we may call the *input/output approach*, relies on operator-theoretic techniques. Among the main contributions to this area, one may cite the foundational work by Zames, Sandberg, Desoer, Safanov, Vidyasagar, and others, starting in the 1960s. In this approach, a "system" is a causal operator $F$ between spaces of signals, and "stability" is taken to mean that $F$ maps bounded inputs into bounded outputs, or finite-energy inputs into finite-energy outputs. More stringent typical requirements in this context are that the gain of $F$ be finite (in more classical mathematical terms, that the operator be bounded), or that it have finite incremental gain (mathematically, that it be globally Lipschitz). The input/output approach has been extremely successful in the robustness analysis of linear systems subject to nonlinear feedback and mild nonlinear uncertainties, and in general in the area that revolves around the various versions of the small-gain theorem. Moreover, geometric characterizations of robustness (gap metric and the like) are elegantly carried out in this framework. Finally, i/o stability provides a natural setting in which to study the classification and parameterization of dynamic controllers.

On the other hand, there is the model-based, or state-space approach, to systems and stability, where the basic object is a forced dynamical system, typically described by differential or difference equations. In this approach, which is consistent with the way we have treated systems in the first part of this paper, there is a standard notion of stability, namely Lyapunov asymptotic stability of the unforced system (GAS). Associated to such a system, there is an operator $F$ mapping inputs (forcing functions) into state trajectories (or into outputs, if partial measurements on states are of interest). It becomes of interest then to ask to what extent Lyapunov-like stability notions for a state-space system are related to the stability, in the senses discussed in the previous paragraph, of the associated operator $F$. It is well-known that, in contrast to the case of linear systems, where there is –subject to mild technical assumptions– an equivalence between state-space and i/o stability, for nonlinear systems the two types of properties are not so closely related. Even for the very special and comparatively simple case of "feedback linearizable" systems, this relation is far more subtle than it might appear at first sight: if one first linearizes a system and then stabilizes the equivalent linearization, in terms of the original system one does not in general obtain a closed-loop system that is input/output stable in any reasonable sense. (However, it is always possible to make a choice of a —usually different— feedback law that achieves such stability, in the linearizable case as well as for all other stabilizable systems; see [35].)

The work initiated in [35] elucidates the links between the two alternative paradigms of stability, i/o and state-space, through the systematic use of the notion of "input to state stability" (ISS). This notion differs fundamentally from the operator-theoretic ones that have been classically used in control theory, first of all because it takes account of initial states in a manner fully compatible with Lyapunov stability. Second, boundedness (finite gain) is far too strong a requirement for general nonlinear operators, and it must be replaced by "nonlinear gain estimates," in which the norms of output signals are bounded by a nonlinear function of the norms of inputs; the definition of ISS incorporates such gains in a natural way.

The ISS concept has already proven useful in many applications; see e.g. [3, 10, 16, 18, 19, 20, 24, 23, 29, 30, 31, 33, 46]. Moreover, this concept has many equivalent versions, which indicates that it is mathematically natural: there are characterizations in terms of dissipation, robustness margins, and classical Lyapunov-like functions; see e.g. [42, 43].

Informally, the ISS property translates into the statement: "no matter what is the initial state, if the inputs are uniformly small, then the state must eventually be small". It gives no useful bounds in the situation in which inputs $u(\cdot)$ are unbounded but still have in some sense "finite energy". In [37], a new notion, that of *integral* input to state stability, iISS for short, was introduced, to model the statement: "no matter what is the initial state, if *integrals* of the inputs are small, then the state must eventually be small". In this part of the paper, we will briefly survey both the ISS and iISS properties.

**Bounded inputs may give unbounded states**

To appreciate the type of problem that one may encounter, consider the following issue. Suppose that (1) is GAS, meaning, more precisely, that the trivial solution $x \equiv 0$ of the differential equation

$$\dot{x} = f_0(x) = f(x, 0) \tag{43}$$

is globally asymptotically stable. Then one would like to know if, for solutions of (1) associated to *non*zero controls, it holds that

$$u(\cdot) \searrow 0 \;\;\Rightarrow\;\; x(\cdot) \searrow 0$$

(the "converging input converging state" property) or that

$$u(\cdot) \text{ bounded} \;\;\Rightarrow\;\; x(\cdot) \text{ bounded}$$

(the "bounded input bounded state" property). Of course, for linear systems $\dot{x} = Ax + Bu$ these implications are always true. Not only that, but one has explicit estimates

$$|x(t)| \;\leq\; \beta(t)\,|x^0| + \gamma \|u_t\|_\infty$$

where
$$\beta(t) = \|e^{tA}\| \to 0 \text{ and } \gamma = \|B\| \int_0^\infty \|e^{sA}\| ds$$

for any Hurwitz matrix $A$, where $u_t$ is the restriction of $u$ to $[0,t]$, thought of as a measurable bounded function which is zero for $s > t$. From these estimates both properties can be easily deduced, see e.g. [36].

These implications fail in general for nonlinear systems, however, as has been often pointed out in the literature (see for instance [47]). As a trivial illustration, take the system

$$\dot{x} = -x + (x^2 + 1)u \tag{44}$$

and the control $u(t) = (2t+2)^{-1/2}$. With $x_0 = \sqrt{2}$ there results the unbounded trajectory $x(t) = (2t+2)^{1/2}$. This is in spite of the fact that the unforced system is GAS. Thus, the converging input converging state property does not hold. Even worse, the bounded input $u \equiv 1$ results in a finite-time explosion. This example is not artificial, as it arises from the simplest case of feedback linearization design. Indeed, given the system

$$\dot{x} = x + (x^2 + 1)u,$$

the obvious stabilizing control law (obtained by first canceling the nonlinearity and then assigning dynamics $\dot{x} = -x$) is

$$u := \frac{-2x}{x^2 + 1} + v$$

where $v$ is the new external input. In terms of this new control (which might be required in order to meet additional design objectives, or may represent the effect of an input disturbance), the closed-loop system is as in (44), and thus is ill-behaved. Observe, however, that if instead of the obvious law just given one would use:

$$u := \frac{-2x}{x^2 + 1} - x + v,$$

then the closed-loop system becomes instead (calling "$v$" again "$u$"):

$$\dot{x} = -2x - x^3 + (x^2 + 1)u.$$

This is still stable when $u \equiv 0$, but in addition it tolerates perturbations far better, since the term $-x^3$ dominates $u(x^2 + 1)$ for bounded $u$ and large $x$. The behavior with respect to such $u$ is characterized qualitatively by the notion of ISS system, to be discussed below.

### Motivation: Finite-Gain under Coordinate Changes

The notions of ISS and iISS arise naturally when one attempts to generalize linear gain analysis to nonlinear systems. In order to motivate the definitions, we first review the classical case of linear systems

$$\dot{x} = Ax + Bu.$$

There are two main norms which one usually imposes on input functions and state trajectories, namely square norms (to quantify energy) or supremum norms (to study worst-case behavior). This gives rise to the following three possible estimates (the case "$L^\infty \to L^2$" is less interesting, being too restrictive):

$$[L^\infty \to L^\infty] \qquad |x(t)| \;\leq\; c\,|x(0)|\,e^{-\lambda t} \;+\; c \sup_{s\in[0,t]} |u(s)| \quad \text{for all } t \geq 0\,, \qquad (45)$$

$$[L^2 \to L^\infty] \qquad |x(t)| \;\leq\; c\,|x(0)|\,e^{-\lambda t} \;+\; c \int_0^t |u(s)|^2\,ds \quad \text{for all } t \geq 0\,, \qquad (46)$$

$$[L^2 \to L^2] \qquad \int_0^t |x(s)|^2\,ds \;\leq\; c\,|x(0)| \;+\; c \int_0^t |u(s)|^2\,ds \quad \text{for all } t \geq 0\,. \qquad (47)$$

(To be precise: there exist positive constants $c$ and $\lambda$ so that, for each measurable essentially bounded control $u$ and each initial condition, the solution of $\dot{x} = Ax + Bu$ is defined for all $t$ and it satisfies the indicated estimates. Actually, most textbooks omit the initial states, but these are the appropriate formulations when nonzero initial states are taken into account.) It is not difficult to see that all these three properties are equivalent for linear systems, and in fact they amount to simply asking that the matrix $A$ be Hurwitz (all eigenvalues have negative real parts).

One may argue, however, that *notions of stability should be invariant under changes of variables.* Let us see what this requirement leads us to. Suppose that we take an origin-preserving state change of coordinates $x = T(z)$ and an origin-preserving change of input variables $u = S(v)$. That is, $T: \mathbb{R}^n \to \mathbb{R}^n$ and $S: \mathbb{R}^m \to \mathbb{R}^m$ are continuous and invertible, have continuous inverses, and $T(0) = 0$ and $S(0) = 0$. Then, there are two functions $\underline{\alpha}, \overline{\alpha} \in \mathcal{K}_\infty$ so that

$$\underline{\alpha}(|z|) \;\leq\; |T(z)| \;\leq\; \overline{\alpha}(|z|)$$

for all $z \in \mathbb{R}^n$, and, similarly, we can write $|S(v)| \leq \overline{\gamma}(|v|)$ for each $v \in \mathbb{R}^m$, for some $\overline{\gamma} \in \mathcal{K}_\infty$. For any input $u$ and initial state $x^0$, and the corresponding trajectory $x(t) = x(t, x^0, u)$, let us write $x(t) = T(z(t))$ and $u(t) = S(v(t))$ for all $t$, and $z^0 = z(0) = T^{-1}(x^0)$. Let us rewrite each of the estimates in terms of the $z$ and $v$ variables.

Estimate (45) translates into:

$$\underline{\alpha}(|z(t)|) \;\leq\; c\,e^{-\lambda t}\overline{\alpha}(|z^0|) \;+\; c \sup_{s\in[0,t]} \overline{\gamma}(|v(s)|) \quad \text{for all } t \geq 0\,.$$

If we use again "$x$" and "$u$" for state and input variables, and we let $\beta(s,t) := ce^{-\lambda t}\overline{\alpha}(s)$ and $\gamma(s) := c\overline{\gamma}(s)$, we are led to an estimate as follows, where $\beta \in \mathcal{KL}$, $\gamma \in \mathcal{K}_\infty$:

$$\underline{\alpha}\,(|x(t)|) \;\leq\; \beta(|x(0)|,t) \;+\; \sup_{s\in[0,t]} \gamma(|u(s)|)\,. \qquad (48)$$

(This would appear to be less general, because $\beta$ has a very special form. Surprisingly, however, any $\mathcal{KL}$ function $\beta$ can be majorized by a function of the special form $c\,e^{-t}\alpha(s)$, with $\alpha \in \mathcal{K}_\infty$; see [37].) On the other hand, estimate (46) translates into:

$$\underline{\alpha}(|z(t)|) \;\leq\; c\,e^{-\lambda t}\overline{\alpha}(|z^0|) \;+\; c \int_0^t \overline{\gamma}(|v(s)|)^2\,ds \quad \text{for all } t \geq 0\,.$$

With $\beta(s,t) := ce^{-\lambda t}\overline{\alpha}(s)$ and $\gamma(s) := c\overline{\gamma}(s)^2$, we have now an estimate as follows::

$$\underline{\alpha}\,(|x(t)|) \;\leq\; \beta(|x(0)|,t) \;+\; \int_0^t \gamma(|u(s)|)\,ds\,. \qquad (49)$$

Finally, estimate (47) translates into:

$$\int_0^t \underline{\alpha}\left(|z(s)|\right) ds \;\leq\; c\overline{\alpha}(|z^0|) + c\int_0^t \overline{\gamma}(|v(s)|)^2\, ds \quad \text{for all } t \geq 0\,.$$

With $\gamma(s) := c\overline{\gamma}(s)^2$ and $\kappa(s) := c\overline{\alpha}(s)$, invariance under change of variables in (47) leads to an estimate of the form:

$$\int_0^t \underline{\alpha}(|x(s)|)ds \;\leq\; \kappa(|x^0|) + \int_0^t \gamma(|u(s)|)ds\,. \tag{50}$$

We now *define* any system (1) to be **input to state stable (ISS)** if there exist functions $\beta \in \mathcal{KL}$ and $\underline{\alpha}, \gamma \in \mathcal{K}_\infty$ so that, for each initial state $x^0$ and input $u$ (measurable and locally essentially bounded), the solution $x(t) = x(t, x^0, u)$ is defined for all $t \geq 0$, and the estimate (48) holds.

Similarly, we say that a system (1) is **integral input to state stable (iISS)** provided that there exist functions $\beta \in \mathcal{KL}$ and $\underline{\alpha}, \gamma \in \mathcal{K}_\infty$ so that, for each initial state $x^0$ and input $u$ (measurable and locally essentially bounded), the solution $x(t) = x(t, x^0, u)$ is defined for all $t \geq 0$, and the estimate (49) holds.

Regarding the last estimate, that in Equation (50) (what one might call a "nonlinear $H_\infty$ estimate" in analogy to the linear case), we do not give it a name, because, as remarked below, it is the same as ISS.

In summary, we have seen that the notions of ISS and iISS arise most naturally when one starts from standard stability notions for linear systems and applies nonlinear changes of variables. (Analogously, although this is not usually mentioned, Lyapunov's asymptotic stability can be seen as the notion which arises, via change of state variables, from exponential stability.)

## 6.1 ISS

Note that the estimate (48) can be rewritten as follows: there exist $\gamma \in \mathcal{K}_\infty$ and $\beta \in \mathcal{KL}$ such that, for all inputs and initial states,

$$|x(t)| \;\leq\; \beta(|x(0)|, t) + \sup_{s \in [0,t]} \gamma(|u(s)|) \tag{51}$$

or, equivalently (with slightly different comparison functions),

$$|x(t)| \;\leq\; \max\left\{\beta(|x(0)|, t),\, \sup_{s\in[0,t]} \gamma(|u(s)|)\right\}\,. \tag{52}$$

(The equivalences follow from elementary properties of $\mathcal{K}_\infty$ functions; for example, $\underline{\alpha}(r) \leq \beta(s,t) + \gamma(t)$ implies $r \leq \underline{\alpha}^{-1}(\beta(s,t) + \gamma(t)) \leq \underline{\alpha}^{-1}(2\beta(s,t)) + \underline{\alpha}^{-1}(2\gamma(t))$. Note that $\underline{\alpha}^{-1}(2\beta(\cdot,\cdot)) \in \mathcal{KL}$ and $\underline{\alpha}^{-1}(2\gamma(\cdot)) \in \mathcal{K}_\infty$.)

**Asymptotic Gain Characterization of ISS**

The estimate (52) implies that an ISS system must have a well-defined *asymptotic gain*, that is, there is some $\gamma \in \mathcal{K}_\infty$ so that

$$\varlimsup_{t \to +\infty} |x(t, x^0, u)| \;\leq\; \gamma\left(\|u\|_\infty\right)$$

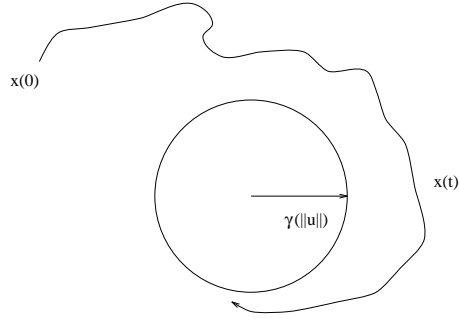

Figure 19: Asymptotic gain: all trajectories tend to ball of radius determined by $\|u\|_\infty$

for all $x^0$ and $u$; see Figure 19. When $u \equiv 0$, (52) reduces to $|x(t)| \leq \beta(|x(0)|, t)$, a GAS estimate, so an ISS system has the property that the zero-input system (43) is GAS. Thus, an ISS system implies the properties: *the zero-input system is GAS* and *there is an asymptotic gain*. The first property is independent of inputs, and the second property involves no estimates on initial states. It is then surprising that ISS is in fact equivalent to the conjunction of these two properties, which represents a **superposition principle** for nonlinear control. Moreover, one can weaken the GAS assumption to simply (neutral) stability of the zero-input system (43), because the asymptotic gain property, when applied to the special case $u \equiv 0$, already implies that all trajectories of (43) converge to zero (attractivity property). The result, which is nontrivial, is as follows, see [43]:

**Theorem H** *The system (1) is ISS if and only if it admits an asymptotic gain and the unforced system (43) is stable.*

We remark that the basic difficulty is in establishing uniform convergence estimates for the states (i.e., constructing the $\beta$ function). If the system were affine in $u$, the proof would follow from a more or less standard weak compactness argument; however, for general systems, the proof is far more delicate.

**Dissipation Characterization of ISS**

Another characterization of ISS is via Lyapunov-like notions. Let us say that a smooth $V : \mathbb{R}^n \to \mathbb{R}$ is an *ISS-Lyapunov function* for the system (1) if $V$ is proper and positive definite, and, for some $\gamma, \alpha \in \mathcal{K}_\infty$,

$$\dot{V}(x,u) \leq -\alpha(|x|) + \gamma(|u|) \quad \forall\, x, u. \tag{53}$$

(Recall that we use the notation "$\dot{V}(x,u)$" to denote $\nabla V(x).f(x,u)$.) In other words, along all trajectories of the system, we have the dissipation inequality

$$V(x(t_2)) - V(x(t_1)) \leq \int_{t_1}^{t_2} w(u(s), x(s))\, ds$$

where, using the terminology of [48, 15], the "supply" function is $w(u,x) = \gamma(|u|) - \alpha(|x|)$. The following fundamental result is from [42]:

**Theorem I** *The system (1) is ISS if and only if it admits an ISS-Lyapunov function.*

It is trivial to see that, if there is an ISS-Lyapunov function, then there is also one which satisfies an estimate of the type (with different comparison functions, of course):

$$\dot{V}(x,u) \leq -\alpha(V(x)) + \gamma(|u|) \tag{54}$$

or even (cf. [31], and also [36], Exercise 5.7.17) one of this form with $\alpha(r) = r$, due to Praly and Wang:

$$\dot{V}(x,u) \leq -V(x) + \gamma(|u|). \tag{55}$$

The proof of one implication in Theorem I is not hard. If we are given a dissipation estimate (54), then either $\alpha(V(x)) \leq 2\gamma(|u|)$ or $\dot{V} \leq -\alpha(V)/2$, and from here one deduces by a comparison theorem that

$$V(x(t)) \leq \max\left\{\beta(V(0),t),\, \alpha^{-1}(2\gamma(\|u\|))\right\}$$

where $\beta(s,t) = y(t)$ and $y$ is obtained by solving the scalar initial value problem

$$\dot{y} = -\frac{1}{2}\alpha(y) + \gamma(u), \quad y(0) = s.$$

An estimate of the type (52) is immediate from here.

The other implication takes a bit more effort. The critical observation is that ISS implies "small gain robustness": if $u(t)$ is for all $t$ sufficiently small compared to $x(t)$, then the term involving $u$ can be dropped from the right-hand side of (52). This leads to the property that there is some $\rho \in \mathcal{K}_\infty$ so that, defining $g(x,d) := f(x,d\rho(x))$, the system $\dot{x} = g(x,d)$ is asymptotically stable uniformly with respect to all measurable $d(\cdot) : \mathbb{R}_{\geq 0} \to \mathbb{R}^m$ with values in a ball of radius one. One may then invoke the converse Lyapunov theorem for robust stability proved in [28] to obtain a smooth $V$ so that

$$\nabla V(x) \cdot g(x,d) \leq -\alpha(x)$$

for all $x$ and $d$, and an ISS-Lyapunov function can be obtained by slightly modifying this $V$.

### ISS Behaves Well Under Series Connections

As a further motivation for the concept of ISS, and as an illustration of the characterizations given, we remark that any cascade (series connection) of ISS systems is again ISS. Consider a cascade connection

$$\begin{aligned} \dot{z} &= f(z,x) \\ \dot{x} &= g(x,u) \end{aligned}$$

as in Figure 20. If each system is separately ISS (the $z$ system having $x$ as its input), then we know that there are ISS-Lyapunov functions for each. Moreover, one can find *matching* such functions $V_i$, in the sense that we have estimates:

$$\begin{aligned} \dot{V}_1(z,x) &\leq \theta(|x|) - \alpha(|z|) \\ \dot{V}_2(x,u) &\leq \tilde{\theta}(|u|) - 2\theta(|x|) \end{aligned}$$

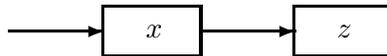

Figure 20: Cascade

(see [41]) and from here one has that $W(x,z) := V_1(z) + V_2(x)$ is an ISS-Lyapunov function for the composite system:

$$\dot{W}(x,z) \leq \tilde{\theta}(|u|) - \theta(|x|) - \alpha(|z|)\,.$$

More generally, if the input $u$ is fed back with "small gain", for instance if $u = k(z)$ is so that $\tilde{\theta}(|u|) \leq (1-\varepsilon)\alpha(|z|)$, then $\dot{W}(x,z) \leq -\theta(|x|) - \varepsilon\alpha(|z|)$ and thus for such feedback the closed-loop system is GAS. This type of reasoning leads to the Small-Gain Theorem for ISS systems due to Jiang, Praly, and Teel, cf. [19], which implies that, under suitable conditions on gains, the closed loop system obtained from an interconnection of two ISS systems $\dot{x} = f(x,z,u)$ and $\dot{z} = g(z,x,v)$, see Figure 21, is itself ISS with respect to $(u,v)$.

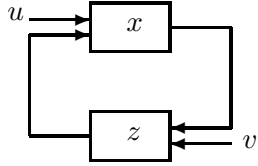

Figure 21: Feedback interconnection

**An Example**

As an example, we look at the problem of angular momentum stabilization of a rigid body controlled by two torques acting along principal axes (for instance, a satellite controlled by two opposing jet pairs). Using $\omega = (\omega_1, \omega_2, \omega_3)$ to denote the angular velocity of a body-attached frame with respect to inertial coordinates, and $I = \mathrm{diag}(I_1, I_2, I_3)$ to denote the principal moments of inertia, we obtain the equations:

$$I\dot{\omega} = \begin{pmatrix} 0 & \omega_3 & -\omega_2 \\ -\omega_3 & 0 & \omega_1 \\ \omega_2 & -\omega_1 & 0 \end{pmatrix} I\omega + \begin{pmatrix} 0 & 0 \\ 1 & 0 \\ 0 & 1 \end{pmatrix} v$$

We assume $I_2 \neq I_3$; then, introducing new state and input coordinates via $(I_2 - I_3)x_1 = I_1\omega_1$, $x_2 = \omega_2$, $x_3 = \omega_3$, $I_2 u_1 = (I_3 - I_1)\omega_1\omega_2 + v_1$, and $I_3 u_2 = (I_1 - I_2)\omega_1\omega_3 + v_2$, we obtain a system on $\mathbb{R}^3$, with controls in $\mathbb{R}^2$:

$$\begin{aligned} \dot{x}_1 &= x_2 x_3 \\ \dot{x}_2 &= u_1 \\ \dot{x}_3 &= u_2\,. \end{aligned}$$

We claim that the following feedback law globally stabilizes the system:

$$u_1 = -x_1 - x_2 - x_2 x_3$$
$$u_2 = -x_3 + x_1^2 + 2x_1 x_2 x_3\,.$$

This feedback was obtained in this way: with $z_2 := x_1 + x_2$, $z_3 := x_3 - x_1^2$, the system becomes:

$$\dot{x}_1 = -x_1^3 + \alpha(x_1, z_2, z_3)$$
$$\dot{z}_2 = -z_2$$
$$\dot{z}_3 = -z_3$$

The $x_1$-subsystem is easily seen to be ISS, because $\deg_{x_1} \alpha \leq 2$ and hence the cubic term dominates for large $x_1$. Thus the cascade, which is the original system that we wanted to stabilize, is also ISS, and, in particular, is GAS when there are no inputs. (In fact, we proved a stronger result, in establishing the ISS property, providing a global robustness result with respect to actuator noise.)

## 6.2 The $H_\infty$-Like Property

Recall that the "$L^2 \to L^2$" stability property in (47) led to an estimate of the form (50), more precisely, for each $x^0$ and $u$ the solution is well-defined for all $t \geq 0$ and (50) holds. (Actually, a weaker notion is obtained if one merely asks that the estimate should hold in the maximal interval of existence.) It is shown in [38] that this property holds *if and only if* the system is ISS. A sketch of a proof is as follows. If the system is ISS, then there is an ISS-Lyapunov function satisfying $\dot{V}(x,u) \leq -V(x) + \gamma(|u|)$ (cf. Equation (55)), so, integrating, we obtain

$$\int_0^t V(x(s))\,ds \;\leq\; \int_0^t V(x(s))\,ds + V(x(t)) \;\leq\; V(x(0)) + \int_0^t \gamma(|u(s)|)\,ds$$

and an estimate of the form (50) holds. Conversely, if such an estimate holds, ons can prove that $\dot{x} = f(x, 0)$ is stable and that an asymptotic gain exists. We show here just the "limit property" $\inf_{t \geq 0} |x(t)| \leq \theta(\|u\|_\infty)$. Indeed, let $\theta := \underline{\alpha}^{-1} \circ \gamma$. Suppose that $\inf_{t \geq 0} |x(t)| > (\underline{\alpha}^{-1} \circ \gamma)(\|u\|)$, so that there is some $\varepsilon > 0$ so that $\underline{\alpha}(x(t)) \geq \varepsilon + \gamma(u(t))$ for all $t \geq 0$. Then, $\int_0^t \underline{\alpha}(x(s))ds \geq \varepsilon t + \int_0^t \gamma(u(s))ds$, which implies $\kappa(|x^0|) > \varepsilon t$ for all $t$, a contradiction.

## 6.3 Integral Input-to-State Stability

The "$L^2 \to L^\infty$" linear stability property (46) led us to the iISS estimate (49). As with the ISS property, there is an elegant dissipation characterization. Let us call a smooth $V : \mathbb{R}^n \to \mathbb{R}$ is an *iISS-Lyapunov function* for the system (1) if $V$ is proper and positive definite, and, for some $\gamma \in \mathcal{K}_\infty$ and some *positive definite* $\alpha$,

$$\dot{V}(x,u) \;\leq\; -\alpha(|x|) + \gamma(|u|) \quad \forall\, x, u \tag{56}$$

(the difference with (53) is that here we do not require $\alpha$ to be class $\mathcal{K}_\infty$). A recent result, obtained in [1] by Angeli, Wang, and the author, is as follows:

**Theorem J** *The system (1) is iISS if and only if it admits an iISS-Lyapunov function.*

Notice that this implies (since a $\mathcal{K}_\infty$ function is in particular positive definite) that every ISS system is also iISS, but the converse is false, because for bounded but large $u$'s one may well have $\dot{V} > 0$. (Intuitively, on the other hand, $\gamma(u) \in \mathcal{L}^1$ means that $\dot{V}$ is "often" negative.) As an example, take the system $\dot{x} = -\tan^{-1} x + u$. This is iISS: take for instance $V(x) = x \tan^{-1} x$, so $\dot{V} \leq -(\tan^{-1} |x|)^2 + 2|u|$, but it is not ISS, since $u(t) \equiv \pi/2$ results in unbounded trajectories.

Several properties are known for the iISS property. For instance, a bilinear system $\dot{x} = (A + \sum_{i=1}^{m} u_i A_i) x + Bu$ is iISS if and only if $A$ is a Hurwitz matrix, see [38]. Many more properties are given in [1]. We pick just one of the most interesting ones here.

A system (1) together with an output function $y = h(x)$ (where $h : \mathbb{R}^n \to \mathbb{R}^p$ is a continuous map, for some $p$, with $h(0) = 0$) will be said here to be *h-dissipative* if there exists a smooth, positive definite, and proper function $V : \mathbb{R}^n \to \mathbb{R}$ such that

$$\dot{V}(x,u) \leq -\alpha(h(x)) + \sigma(|u|)$$

for all $x, u$, and *weakly h-detectable* if, for all trajectories, $y(t) = h(x(t)) \equiv 0$ implies that $x(t) \to 0$ as $t \to \infty$.

**Theorem K** *The system (1) is iISS if and only if there is some output $h$ so that the system is weakly h-detectable and h-dissipative.*

## 7 Other Notions Related to ISS

There are several notions, all of great interest in control, which arise when applying coordinate changes to classical linear concepts. We very briefly mention some of them here.

Consider first systems $\dot{x} = f(x), y = h(x)$ with outputs but no inputs. Dualizing ISS by simply replacing "$u$" by "$y$", there results the following estimate:

$$|x(t, x^0)| \leq \beta(|x^0|, t) + \sup_{s \in [0,t]} \gamma(|y(s)|).$$

A system for which this property holds (along all maximal solutions) is said to be **output to state stable (OSS)**. This property is precisely the appropriate version of detectability with respect to zero states, see [44]. More generally, one could have inputs *and* outputs, considering systems $\dot{x} = f(x,u), y = h(x)$. The obvious notion of detectability is then given by estimates of the type

$$x(t) \leq \beta(|x^0|, t) + \sup_{s \in [0,t]} \gamma_1(|u(s)|) + \sup_{s \in [0,t]} \gamma_2(|y(s)|)$$

("IOSS", for input/output to state stability). Several dissipation characterizations are discussed in [44, 22].

A different notion is that of **input-to-output stability (IOS)**, where we ask that all trajectories of a system with outputs must satisfy

$$h(x(t, x^0, u)) \leq \beta(|x^0|, t) + \sup_{s \in [0,t]} \gamma(|u(s)|)$$

and

$$|x(t, x^0, u)| \leq \max\left\{\sigma(|x^0|),\, \sup_{s \in [0,t]} \sigma(u(s))\right\}$$

for all $t$, $u$, and $x^0$, and for appropriate comparison functions. This is closely related to the concept of *partial stability* in differential equations (special case, when $h$ is a projection). It is the correct notion, it would appear, for formulating regulation and tracking problems. See [45] for Lyapunov-like dissipation characterizations.

**Acknowledgements**

I wish to thank Francesca Albertini, Misha Krichman, Daniel Liberzon, Gene Ryan, and Yuan Wang, for many useful comments and suggestions.

# References

[1] Angeli, D., E.D. Sontag, and Y. Wang, "A characterization of integral input to state stability," submitted. (Summary as "Remarks on integral input to state stability," in *Proc. IEEE Conf. Decision and Control, Tampa, Dec. 1998*, IEEE Publications, 1998, to appear.)

[2] Artstein, Z., "Stabilization with relaxed controls," *Nonlinear Analysis, Theory, Methods & Applications* **7**(1983): 1163-1173.

[3] Battilotti, S., "Robust stabilization of nonlinear systems with pointwise norm bounded uncertainties: A control Lyapunov approach," *IEEE Trans. Autom. Control*, to appear.

[4] Brockett, R.W., "Asymptotic stability and feedback stabilization," in *Differential Geometric Control Theory* (R.W. Brockett, R.S. Millman, and H.J. Sussmann, eds.), Birkhauser, Boston, 1983, pp. 181-191.

[5] Clarke, F.H., Yu.S. Ledyaev, and R.J. Stern, "Asymptotic stability of differential inclusions," submitted.

[6] Clarke, F.H., Yu.S. Ledyaev, R.J. Stern, and P. Wolenski, *Nonsmooth Analysis and Control Theory*, Springer-Verlag, New York, 1998.

[7] Clarke, F.H., Yu. S. Ledyaev, E.D. Sontag, and A.I. Subbotin, "Asymptotic controllability implies feedback stabilization," *IEEE Trans. Automat. Control* **42**(1997): 1394-1407.

[8] Coron, J-M., "Global asymptotic stabilization for controllable systems without drift," *Math of Control, Signals, and Systems* **5**(1992): 295-312.

[9] Coron, J.-M., and L. Rosier, "A relation between continuous time-varying and discontinuous feedback stabilization," *J.Math. Systems, Estimation, and Control* **4**(1994): 67-84.

[10] Christofides, P.D., and A.R. Teel, "Singular perturbations and input-to-state stability," *IEEE Trans. Automat. Control* **41**(1996): 1645-1650.

[11] Freeman, R., and P. V. Kokotović, *Robust Nonlinear Control Design*, Birkhäuser, Boston, 1996.

[12] Filippov, A.F., *Differential Equations with Discontinuous Right-Hand Side*, Nauka, Moscow, 1985 (in Russian). and Kluwer, Dordrecht, 1988 (English translation).

[13] Hájek, O., "Discontinuous differential equations. I, II", *J. Diff. Equations* **32** (1979), 149-170, 171-185.

[14] Hermes, H., "Discontinuous vector fields and feedback control," in *Differential Equations and Dynamical Systems*, Academic Press, New York, 1967, pp. 155-165.

[15] Hill, D.J., "Dissipative nonlinear systems: Basic properties and stability analysis," *Proc. 31st IEEE Conf. Dec. and Control*, Tucson, Arizona, 1992, pp. 3259-3264.

[16] Hu, X.M., "On state observers for nonlinear systems," *Systems & Control Letters* **17** (1991), pp. 645–473.

[17] Isidori, A., *Nonlinear Control Systems, Third Edition*, Springer-Verlag, London, 1995.

[18] Isidori, A., "Global almost disturbance decoupling with stability for non minimum-phase single-input single-output nonlinear systems," *Systems & Control Letters* **28** (1996), pp. 115–122.

[19] Jiang, Z.-P., A. Teel, and L. Praly, "Small-gain theorem for ISS systems and applications," *Mathematics of Control, Signals, and Systems,* **7** (1994), pp. 95–120.

[20] Khalil, H.K., *Nonlinear Systems*, Prentice-Hall, Upper Saddle River, NJ, second ed., 1996.

[21] Krasovskii, N.N., and A.I. Subbotin, *Positional differential games*, Nauka, Moscow, 1974 [in Russian]. French translation *Jeux differentiels*, Editions Mir, Moscou, 1979. Revised English translation *Game-Theoretical Control Problems*, Springer-Verlag, New York, 1988.

[22] Krichman, M., E.D. Sontag, and Y. Wang, "Input/output to state stability," in preparation.

[23] Krstić, M., I. Kanellakopoulos, and P. V. Kokotović, *Nonlinear and Adaptive Control Design*, John Wiley & Sons, New York, 1995.

[24] Krstić, M. and H. Deng, *Stabilization of Uncertain Nonlinear Systems*, Springer-Verlag, London, 1998.

[25] Lafferriere, G.A., and E.D. Sontag, "Remarks on control lyapunov functions for discontinuous stabilizing feedback," *Proc. IEEE Conf. Decision and Control, San Antonio*, IEEE Publications, 1993, pp. 306-308.

[26] Ledyaev, Yu.S., and E.D. Sontag, "A remark on robust stabilization of general asymptotically controllable systems," in *Proc. Conf. on Information Sciences and Systems (CISS 97)*, Johns Hopkins, Baltimore, MD, March 1997, pp. 246-251.

[27] Ledyaev, Yu.S., and E.D. Sontag, "A Lyapunov characterization of robust stabilization," *J. Nonlinear Analysis*, to appear. Summarized verion in "Stabilization under measurement noise: Lyapunov characterization," *Proc. American Control Conf.*, Philadelphia, June 1998, pp. 1658-1662.

[28] Lin, Y., E.D. Sontag, and Y. Wang, "A smooth converse Lyapunov theorem for robust stability," *SIAM J. Control and Optimization* **34**(1996): 124-160.

[29] Lu, W.M., A class of globally stabilizing controllers for nonlinear systems, *Systems & Control Letters* **25** (1995), pp. 13–19.

[30] Lu, W..M., "A state-space approach to parameterization of stabilizing controllers for nonlinear systems," *IEEE Transactions on Automatic Control* **40** (1995): 1576–1588.

[31] Praly, L., and Y. Wang, "Stabilization in spite of matched unmodelled dynamics and an equivalent definition of input-to-state stability," *Math. of Control, Signals, and Systems* **9** (1996): 1-33.

[32] Ryan, E.P., "On Brockett's condition for smooth stabilizability and its necessity in a context of nonsmooth feedback," *SIAM J. Control Optim.* **32**(1994): 1597-1604.

[33] Sepulchre, R., M. Jankovic, and P.V. Kokotovic, "Integrator forwarding: a new recursive nonlinear robust design," *Automatica* **33**(1997): 979-984.

[34] Sontag E.D., "A Lyapunov-like characterization of asymptotic controllability," *SIAM J. Control and Opt.* **21**(1983): 462-471.

[35] Sontag, E.D., 'Smooth stabilization implies coprime factorization", *IEEE Trans. Automat. Contr.* **34**(1989): 435-443.

[36] Sontag, E.D., *Mathematical Control Theory, Deterministic Finite Dimensional Systems*, Second Edition, Springer-Verlag, New York, 1998.

[37] Sontag, E.D,. "Comments on integral variants of ISS," *Systems & Control Letters* **34**(1998): 93-100.

[38] Sontag, E.D., "Comments on integral variants of ISS," *Systems and Control Letters* **34**(1998): 93-100.

[39] Sontag, E.D., and H.J. Sussmann, "Remarks on continuous feedback," in *Proc. IEEE Conf. Decision and Control, Albuquerque, Dec. 1980*, IEEE Publications, Piscataway, pp. 916-921.

[40] Sontag, E.D., and H.J. Sussmann, "Nonsmooth Control Lyapunov Functions," in *Proc. IEEE Conf. Decision and Control, New Orleans*, IEEE Publications, 1995, pp. 2799-2805.

[41] Sontag, E.D., and A. Teel, "Changing supply functions in input/state stable systems," *IEEE Trans. Autom. Control* **40**(1995): 1476-1478.

[42] Sontag, E.D., and Y. Wang, "On characterizations of the input-to-state stability property," *Systems & Control Letters* **24** (1995): 351-359.

[43] Sontag, E.D., and Y. Wang, "New characterizations of the input to state stability property," *IEEE Trans. on Autom. Control* **41** (1996): 1283-1294.

[44] Sontag, E.D., and Y. Wang, "Output-to-state stability and detectability of nonlinear systems," *Systems and Control Letters* **29**(1997): 279-290.

[45] Sontag, E.D., and Y. Wang, "A notion of input to output stability," *Proc. European Control Conf.*, Brussels, July 1997, Paper WE-E A2, CD-ROM file ECC958.pdf, 6 pages. (Journal version in preparation.)

[46] Tsinias, J., "Input to state stability properties of nonlinear systems and applications to bounded feedback stabilization using saturation," *ESAIM Control Optim. Calc. Var.* **2**(1997): 57-85.

[47] Varaiya, P.P. and R. Liu, "Bounded-input bounded-output stability of nonlinear time-varying differential systems," *SIAM J.Control* **4**(1966): 698-704.

[48] Willems, J.C., "Mechanisms for the stability and instability in feedback systems," *Proc. IEEE* **64** (1976): 24-35.